\documentclass[12pt]{amsart}

\usepackage{amsmath}
\usepackage{amsxtra}
\usepackage{amscd}
\usepackage{amsthm}
\usepackage{amsfonts}
\usepackage{amssymb}
\usepackage[utf8]{inputenc}
\usepackage[french]{babel}
\usepackage{bbm}
\usepackage{hyperref}
\usepackage[mathscr]{eucal}

\numberwithin{equation}{section}

\newtheorem{lemme}[equation]{Lemme}
\newtheorem{proposition}[equation]{Proposition}

\theoremstyle{remark}

\newcommand\hfld[2]{\ \smash{\mathop{\hbox to 7mm{\rightarrowfill}}
     \limits^{\scriptstyle#1}_{\scriptstyle#2}}\ }
\newcommand\hflg[2]{\smash{\mathop{\hbox to 7mm{\leftarrowfill}}
     \limits^{\scriptstyle#1}_{\scriptstyle#2}}}

\newcommand\ogg{\leavevmode\raise.3ex\hbox{$\scriptscriptstyle\langle\!\langle$}\,}
\newcommand\fgg{\leavevmode\raise.3ex\hbox{$\scriptscriptstyle\,\rangle\!\rangle$}}

\newcommand\Oc{{\mathscr O}}

\newcommand\calF{{\mathcal F}}

\newcommand\calO{{\mathcal O}}

\newcommand\calG{{\mathscr G}}

\newcommand\mathfrakg{{\mathfrak g}}

\newcommand\mathfrakp{{\mathfrak p}}

\newcommand\Af{{\mathfrak A}}

\newcommand\build[3]{\mathrel{\mathop{\kern
0pt#1}\limits_{\textstyle #2}^{\textstyle #3}}}

\newcommand\CC{{\mathbbm C}}

\newcommand\FF{{\mathbbm F}}

\newcommand\ZZ{{\mathbbm Z}}
\newcommand\NN{{\mathbbm N}}
\newcommand\QQ{{\mathbbm Q}}

\newcommand\bbP{{\mathbbm P}}
\newcommand\bbA{{\mathbb A}}

\newcommand\bbK{{\mathbbm K}}
\newcommand\bbL{{\mathbbm L}}

\newcommand\rmH{{\rm H}}

\newcommand\GL{{\rm GL}}

\newcommand\Gal{{\rm Gal}}

\newcommand\pr{{\rm pr}}

\newcommand\tr{{\rm tr}}

\newcommand\qd{{\rm qd}}

\newcommand\Spec{{\rm Spec}}

\newcommand\Gm{{\mathbb G}_{m}}
\newcommand\SL{{\rm SL}}

\newcommand\Bun{{\rm Bun}}
\newcommand\Hck{{\rm Hecke}}

\newcommand\hookr\hookrightarrow
\newcommand\hookl\hookleftarrow

\newcommand\isom{\,\smash{\mathop{\hbox to 5mm{\rightarrowfill}}
    \limits^\sim}\,}

\newcommand\rs{{\rm rs}}

\newcommand\reg{{\rm reg}}

\newcommand\Ql{{\overline\QQ_\ell}}

\newcommand\sep{\rm sep}
\newcommand\nr{\rm sep}
\newcommand\on{\operatorname}
\renewcommand\Re{\operatorname{Re}}
\newcommand\N{\on{N}}
\newcommand\st{\on{st}}
\newcommand\der{\on{der}}
\newcommand\res{\on{res}}
\newcommand\meas{\on{meas}}
\newcommand\geom{\on{geom}}
\newcommand\norm{\on{norm}}
\newcommand\spl{\on{sp}}

\begin{document}
	\title[Formule des traces]
	{Formule des traces et fonctorialité: \\
	le début d'un programme}
	\author{Edward Frenkel, Robert Langlands et Ngô Bao
          Ch\^au}\thanks{La recherche de Edward Frenkel
a été subventionnée par DARPA (grant n. HR0011-09-1-0015) et Fondation
Sciences mathématiques de Paris. \\ La recherche de Ngô Bao Ch\^au
a été subventionnée par Simonyi Foundation.}
	\maketitle
	
	\section*{Introduction}	
	
	L'un de nous, Langlands,
encouragé par les travaux d'un deuxième, Ngo, sur le lemme fondamental dont l'absence d'une 
démonstration pendant plus de deux décennies entravait à maints égards tout progrès sérieux de 
la théorie analytique des formes automorphes, avait esquissé un programme pour établir la 
fonctorialité, l'un des deux objectifs principaux de cette théorie. Le troisième, 
Frenkel,
a 
observé que quelques 
idées et formules extraites de la forme géométrique de la correspondance (parfois dite réciprocité 
ou correspondance de Langlands) prévue entre formes automorphes et représen\-tations galoisiennes 
appuient fortement la stratégie envisagée. Ce sont là les trois points de départ de cet article.

Le lemme fondamental est maintenant acquis grâce aux résultats de \cite{N} et par 
conséquent la formule des traces 
stable à portée de main. Quoique conscients des obstacles qui restent et de leur difficulté et 
sans être encore en état de les surmonter, nous entreprenons dans cet article et ceux qui le 
suivront la première étape du programme esquissé dans \cite{L4} qui a pour but d'établir les 
liens envisagés entre les formes automorphes et la géométrie algébrique, diophantienne ou non. 
Nous sommes conscients que la première étape, pour ne pas parler de celles qui suivent, n'est pas 
tenue pour réaliste par la plupart (presque tous) des arithméticiens et même des spécialistes 
de la théorie des formes automorphes. Nous espérons néanmoins qu'ils lisent
attentivement, avec un esprit ouvert, cet article
et ceux 
qui le suivront.

Quoique dans ce premier article nous n'envisageons que la théorie sur les corps de nombres et les 
corps de fonctions de dimension un sur un corps fini, nous voulons dans la suite traiter, non pas 
d'une façon uniforme mais d'une façon parallèle, la théorie pour les courbes
algébriques sur le corps des nombres complexes, mais il y a toujours des points obscurs
qu'il reste à éclaircir.

Nous pouvons expliquer en quelques mots le contenu de cet article. Dans la première section
nous résumons brièvement les idées de \cite{L3} mais en utilisant en plus 
des dérivées logarithmiques des fonctions $L$ automorphes ces fonctions
elles-mêmes. Nous observons que cette façon d'aborder la fonctorialité a des relents
de la théorie classique des corps de classes. Dans la deuxième section, qui est la seule
dans laquelle le corps $F$ est restreint à un corps de fonctions, nous vérifions
un lemme géométrique qui se justifie en anticipant
la correspondance entre formes
automorphes pour le corps $F$ de fonctions sur une courbe sur $\FF_q$
et représentations $\ell$-adiques de $\Gal(F^{\rm sep}/F)$. Ensuite, après avoir rappelé
et reformulé la formule des traces stable dans la troisième
section, nous introduisons dans la quatrième ce que nous appelons l'adélisation de la formule
des traces. Dans la formule des traces interviennent des mesures. Dans le passé 
ces mesures ont rendu difficile l'utilisation efficace
de la formule car elles introduisent des facteurs qui
empêchent l'exploitation de la méthode de Poisson, à peu près 
la seule méthode qui semble prometteuse. Il y a pourtant une observation
importante:
avec des hypothèses tout à fait anodines, la partie la plus 
intéressante de la forme stabilisée de la formule, 
à savoir la somme sur les éléments réguliers, est à peu près une somme sur un 
espace vectoriel $V$ de dimension finie sur $F$ d'une
fonction sur $V\otimes\bbA_F$ dont le comportement
permet l'application de la formule de Poisson
pour la paire $V\subset V\otimes\bbA_F$. Il nous faut utiliser cependant
des travaux antérieurs de Kottwitz sur 
les facteurs donnés par les mesures, car ses formules permettent de conclure que ces facteurs
sont constants, ce qui est à notre avis un miracle. Le lecteur est invité
à se pencher sur les renvois aux travaux de Kottwitz et sur ce que 
nous en déduisons et d'y réfléchir. Nous soulignons que ces facteurs ne se
simplifient que
pour la formule stable. Nous ajoutons que pour utiliser la formule
de Poisson il nous faudra tronquer l'espace $V\otimes\bbA_F$ et
à cette fin nous avons emprunté une idée de Jayce Getz, une idée qu'il
utilise lui-même dans un autre contexte mais toujours dans le cadre
de la formule des traces. Nous lui
sommes reconnaissants de nous l'avoir expliquée. 

Il est impossible de souligner toute l'importance de cette possibilité
de traiter les problèmes analytiques que pose la formule des traces.
Jusqu'à présent personne ne s'en est aperçu.
Nous n'abordons pas cependant dans cet article l'analyse des sommes de Poisson,
mais nous vérifions dans la dernière section que le premier terme, $\hat\theta(0)$,
de la somme discrète sur $F$ de la transformée $\hat\theta$ qui y apparaît est la contribution
dominante, celle des représentations automorphes de dimension un. C'est un très bon
signe. Nous ajoutons cependant que cette contribution, bien que dominante, n'est guère la
plus importante.

L'un de nous est particulièrement content que cet article apparaisse dans un numéro dédié à Paulo 
Ribenboim, car c'est à la suite d'une demande de sa part que les premiers balbutiements de l'endoscopie ont 
été rédigés et ont paru dans le {\it Journal canadien de mathématiques}.

\section{Pôles des fonctions $L$ et fonctorialité}\label{L} 

Pour les définitions de base des fonctions $L$ rattachées aux représent\-ations automorphes nous 
renvoyons au livre {\it Automorphic forms, representations
and $L$-functions} \cite{BC}. Il y a d'excellentes introductions plus récentes à la théorie des 
représentations automorphes, telles que \cite{AEK}, qui toutefois ne traitent pas les fonctions $L$. 
Le lecteur aura besoin en lisant cet article de quelque compréhension de la formule des traces 
stable et des applications prévues, donc des formes automorphes au-delà de l'endoscopie. 
À certains égards et malgré les contributions héroïques et fondamentales de Arthur, la formule des 
traces stable n'existe toujours que sous une forme rudimentaire \cite{A2,A3,K2,L2}. Nous voulons expliquer
dans cet article et ceux qui le suivent ce qu peuvent être ses objectifs
et comment avec son aide on peut espérer les réaliser. Nous n'avançons qu'à tâtons mais pour 
orienter le lecteur nous commençons avec une description du fond.  

Les fonctions $L$ sont des produits (pris sur toutes les places $v$ finie ou infinie)
$$
L(s,\pi,\rho)=\prod_v L(s,\pi_v,\rho)
$$ 
rattachés à des formes automorphes ou plutôt à des $L$-paquets de formes automorphes d'un groupe 
$G$ et à une représentation $\rho$ de son $L$-groupe.

Arthur (\cite{A1}) a proposé une classification des formes automorphes qui ne sera réalisé qu'à 
partir de la fonctorialité, qui elle-même ne peut sans doute être établie qu'en même temps que
la classification. Pour comprendre 
notre stratégie il faut comprendre au moins quels seront les éléments d'une telle classification. 
Pour l'établir en général, on s'attend à utiliser la formule des traces et des récurrences. Nous 
verrons dans les prochaines pages comment tenir compte des conséquences de cette classification en 
maniant cette formule.

On s'attend à pouvoir rattacher à une  représentation automorphe $\pi=\pi_G$ plusieurs objets dont 
d'abord un homomorphisme $\phi=\phi_\pi$ de $\SL(2)$ dans ${}^LG$, ensuite pour presque toute 
place $v$ une classe de conjugaison $\{A_G(\pi_v)\}=\{A(\pi_v)\}$ dans le centralisateur connexe 
${}^\lambda G_{\phi,v}$ de $\phi(\SL(2))$ dans ${}^LG_v$ et même pour toute place $v$ un 
homomorphisme du groupe de Weil local dans ${}^\lambda G_{\phi,v}$. On s'attend même à pouvoir 
construire un groupe de Galois automorphe et un homomorphisme $\xi$ de ce dernier dans le 
centralisateur global ${}^\lambda G_{\phi}$ qui engendre les classes locales $\{A(\pi_v)\}$. Il 
s'agit d'un projet mais d'un projet qui est censé mener à sa propre réalisation.

En fait ce groupe de Galois serait inutilement gros et guindé.  
Plus utiles seraient les groupes ${}^\lambda H_\pi$ définis par les clôtures de l'image de
l'application hypothétique $\xi$ pour la 
topologie de Zariski. Ils sont plus primordiaux que les homomorphismes $\xi$ et 
on s'attend, le cas échéant, à pouvoir définir un groupe de Galois automorphe après avoir déterminé les ${}^\lambda 
H_\pi$. Le problème de la détermination de ${}^\lambda H_\pi$ fut entamé dans \cite{L3}. On cherche 
même un groupe réductif $H$ sur $F$, un homomorphisme surjectif 
$$
\psi:\,{}^LH\rightarrow {}^\lambda H\subset {}^\lambda G_\phi
$$ 
à noyau central, et une représentation $\pi_H$ telle que 
$$
\{A(\pi_v)\}=\{\psi(A(\pi_{H,v})\}
$$ 
pour presque tout $v$. Selon la classification proposée par Arthur, il y a une classe particulière 
de représentations, celles, dites de type Ramanujan, qui satisfont à la conclusion de la 
conjecture de Ramanujan. On veut que $\pi_H$ soit de type Ramanujan, c'est-à-dire que les classes 
$\{A(\pi_{H,v})\}$ soient unitaires. Le problème abordé dans \cite{L3} fut la construction, dans le 
cadre de la formule des traces, des données $\phi$, $H$, $\pi_H$ et $\psi$ à partir de $\pi$. Nous 
remarquons en passant qu'on ne s'attend pas à ce que ces données soient uniques. Il y aura des 
questions de multiplicité liées à ces constructions qui ont été abordées d'une façon concrète par 
Song Wang \cite{WS}.

On s'attend à pouvoir démontrer deux choses en utilisant une récur\-rence convenable avec la formule 
des traces stable et en se fondant sur le principe de fonctorialité : d'une part, que les classes 
$\{A(\pi_{H,v})\}$ sont unitaires de sorte que les valeurs propres de $\rho_H(A(\pi_{H,v}))$ sont 
de valeur absolue $1$ pour toute représentation $\rho_H$ de ${}^LH$ et d'autre part, que les fonctions 
$L(s,\pi_H,\rho_H)$ sont holomorphes pour $\Re s>1$ avec au plus un nombre fini de pôles sur la 
droite $\Re(s)=1$ et aucun zéro sur cette droite critique. (Si le corps $F$ est un corps de 
fonctions sur un corps fini à $q$ éléments, on compte les pôles modulo $2\pi i/\ln q$.) Observons 
que l'influence d'un nombre fini de places dans le produit eulérien est très faible de sorte que le comportement décrit est 
au fond celui des fonctions $L$ partielles 
\begin{equation}\label{1.1}
L_S(s,\pi,\rho)=\prod_{v\notin S}L(s,\pi_{v},\rho).
\end{equation}

La représentation $\rho$ et l'homomorphisme $\phi\times\psi$ définissent une repré\-sentation du 
produit $\SL(2)\times{}^LH$, laquelle se décompose en somme directe,
\begin{equation}\label{RL.1.0}
\bigoplus_j \sigma_j\otimes\rho^j_H,
\end{equation}
où chaque $\sigma_j$ est irréductible. Supposons que la représentation $\sigma_j$ de $\SL(2)$ soit 
de dimension $m_j+1$ et par conséquent que son poids maximal soit $\mu_j:(1,-1)\rightarrow m_j$. 
La fonction $L$ se décompose en un produit, à savoir
\begin{equation}\label{1.2}
L_S(s,\pi,\rho)=\prod_j\prod_i L_{S}(s+i,\pi_H,\rho^j_H)
\end{equation}
où 
$$
i\in \left\{\frac{m_j}{2},\frac{m_j}{2}-1,\ldots,-\frac{m_j}{2}\right\}.
$$
Évidemment le comportement de ces fonctions dépend fortement des entiers $m_j$. Pour $m_j=0$ la 
droite critique est $s=1$, la bande critique est $0\leq\Re s\leq 1$, et son centre $\Re s=1/2$. 
Pour $m_j>0$, les singularités commencent plus tôt pour $s$ décroissant. Si toutefois $m_j>0$, 
alors $\dim H<\dim G$ de sorte que nous pouvons supposer par récurrence, soit sur $\dim G$, soit 
sur $\dim\rho$, que nous comprenons le comportement des fonctions $L_S(s,\pi_H,\rho^j_{H})$
et que nous avons démontré les hypothèses de Arthur pour $H$, en particulier l'hypothèse
de Ramanujan. Nous rappelons que cette dernière hypothèse équivaut à l'hypothèse que,
pour les représentations telles que $\phi$ est trivial, $L(s,\pi,\rho)$
se prolonge jusqu'à la droite $\Re s=1$. 

Pour être plus précis, de la formule des traces on déduit d'abord une formule pour les
sommes
\begin{equation}\label{RL.1.1}
\sum_\pi\prod_{v\in S}\tr(\pi_v(f_v))L_S(s,\pi,\rho).
\end{equation} 
Les fonctions $f_v$, lisses et à support compact, sont à peu près arbitraires
et nous permettent d'isoler à la fin de l'argument les répresentations $\pi$,
ou au moins leurs classes stables,
les séparant les unes des autres. En plus il est implicite dans de telles
expressions que la somme se fait sur les $\pi$ non ramifiés en-dehors de $S$.
Puisque l'on connaît les sous-groupes ${}^LH$ et les $\phi$
tels que l'image de $\phi$ commute avec ${}^LH$, nous saurons par récurrence déduire de la formule
des traces pour les groupes $H$ une formule pour les contributions à \eqref{RL.1.1}
qui proviennent d'un $H$ qui n'est pas $G$ lui-même. Il faudra toutefois tenir compte de la
possibilité qu'une seule $\pi$ provient de plusieurs $H$, par exemple de $H_1$
et $H_2$ tels que ${}^L H_1\subset {}^LH_2$. Des exemples suggèrent aussi que
$\pi$ puisse provenir de deux sous-groupes isomorphes mais non conjugués,
un phénomène qui semble être lié à l'existence de multiplicités plus grandes que $1$ \cite{WS}.
On attend des éclaircissements au fur et à mesure que les recherches
progressent. Tout $\pi$
qui apparaît dans la différence sera alors hypothétiquement de type Ramanujan de sorte que la
différence est censée être une fonction de $s$ holomorphe dans le domaine $\Re s>1$ et le problème
principal sera de le montrer et d'en déduire la conjecture de Ramanujan. Évidemment ce problème 
ne sera guère facile. Nous l'abordons
d'un côté facile dans la partie \S5. Le résultat n'est pas dépourvu d'intérêt! Il est évident
qu'une démonstration de la fonctorialité en général est implicite dans ces propos.

Il est préférable d'écarter les $H$ différents de $G$ en
deux étapes. On soustrait d'abord les contributions des paires $(\phi,\psi)$
telles que $\phi$ n'est pas trivial. Pour elles les pôles des fonctions \eqref{1.2}
apparaissent pour $\Re s$ plus grand de sorte que les comparaisons, qui se font
dans des petits intervalles juste à la droite des pôles, puissent se faire
successivement et indépendamment. L'objectif est d'en tenir compte à tour de rôle
en passant de droite à gauche. À la fin on aura isolé les $\pi$ pour lesquels
$\phi$ est censé être trivial et démontré que la fonctorialité prédite par l'hypothèse
de Arthur est au moins vraie pour $\phi$ non trivial. Les représentations $\pi$ pour lesquelles $\phi$
est trivial sont celles pour lesquelles il faut démontrer la conjecture de Ramanujan,
donc pour lesquelles il est nécessaire de démontrer que pour tout $\rho$
la fonction $L(s,\pi,\rho)$ se prolonge
à la région $\Re s>1$.

De toute façon nous aurons une somme semblable à \eqref{RL.1.1} sauf que dans la somme 
n'apparaîtront que les représentations $\pi$ qui sont censées être de type
Ramanujan.
Donc de cette façon nous aurons isolé les représentations de $G$ de type Ramanujan.
Il faudra ensuite, pour un 
$\rho$ donné, isoler celles dont la fonction $L(s,\pi,\rho)$ a un pôle 
d'un ordre donné en $s=1$, 
ou en n'importe quel autre point donné sur la droite critique $\Re s=1$. Pour isoler une
représentation donnée il faudra profiter du choix arbitraire des fonctions $f_v$
dans \eqref{RL.1.1} de l'ensemble $S$. 

Nous pouvons, à 
partir de la formule des traces, exprimer des sommes de $\tr(\pi(f))$ sur toutes les 
représentations automorphes $\pi$ pour des fonctions $f=\prod_vf_v$ à peu près arbitraires.  Au 
début (\cite{L3}) il semblait préférable de mettre dans \eqref{RL.1.1} non pas les fonctions $L$ de \eqref{1.1} qui 
sont des produits mais leurs dérivées logarithmiques, car pour ces dérivées les résidus sont additifs. 
Cela mène à des sommes semblables à la somme 
$$
\sum_{p\leq X}\ln p,
$$
rencontrée dans la démonstration du théorème des nombres premiers. Nous les appellerons des sommes 
arithmétiques. Bien que nous passerons plus tard aux sommes géométriques
\eqref{RL.1.1}, nous commençons avec les sommes 
arith\-mé\-tiques.

Si $\pi$ est non ramifié en dehors de $S$, la dérivée logarithmique de \eqref{1.2} prise avec un 
signe négatif est égale à
\begin{equation}\label{1.3}
-\frac{L'_S(s,\pi,\rho)}{L_S(s,\pi,\rho)}=
-\sum_{v\notin S}\frac{L'_v(s,\pi,\rho)}{L_v(s,\pi,\rho)}=
\sum_{v\notin S}\ln\N\mathfrakp_v\sum_n\frac{\tr\rho^n(A(\pi_v))}{N\mathfrakp_v^{ns}}.
\end{equation}
Si l'on avait la conjecture de la conjecture de
Ramanujan en main, on pourrait écrire \eqref{1.3} comme une somme,
\begin{equation}\label{1.4}
\sum_v\ln\N\mathfrakp_v\frac{\tr\rho(A(\pi_v))}{N\mathfrakp_v^{s}}
+\sum_{v\notin S}O(\frac{1}{N\mathfrakp_v^2}).
\end{equation}
Cependant ce résultat n'est pas disponible dès notre point de départ. C'est un objectif.

Comme déjà expliqué, on ne s'attend pas  à calculer \eqref{1.3} pour une seule 
représentation $\pi$. Ce que la formule des 
traces stable donne est
\begin{equation}\label{1.5}
\sum_{v\notin S}\ln\N\mathfrakp_v\sum_n\left\{
\sum_{\pi^{\st}}m(\pi^{\st})\prod_{v\in S}\tr\pi_v^{\st}(f_v)
\frac{\tr\rho^n(A(\pi_v))}{N\mathfrakp_v^{ns}}\right\},
\end{equation}
ou
\begin{equation}\label{1.6}
\sum_{v\notin S}\ln\N\mathfrakp_v\left\{
\sum_{\pi^{\st}}m(\pi^{\st})\prod_{v\in S}\tr\pi_v^{\st}(f_v)
\frac{\tr\rho(A(\pi_v))}{N\mathfrakp_v^{s}}\right\},
\end{equation}
car on peut trouver dans l'algèbre de Hecke en la place $v$ une fonction $K^{n}_v$ telle que
$$
\tr\rho^n(A(\pi_v))=\tr\pi_v(K^n_v).
$$
Les fonctions $f_v$ sur $G(F_v)$ sont lisses et à support compact mais arbitraires
à part ces deux conditions. Cela nous permet  de séparer à certaines fins les représentations qui interviennent, 
ou  plutôt les $L$-paquets , notés $\pi^{\st}$, car ce sont eux qui interviennent dans la formule 
des traces stable. La multiplicité $m(\pi^{\st})$ est une multiplicité stable, une notion dont la 
définition générale et précise ne sera sans doute donnée qu'au fur et à mesure de la création 
d'une formule des traces stable. Pour le moment, elle n'a pas été donnée que dans des cas 
particuliers; voir par exemple \cite{LL}. Dans tous les cas, elle sera un nombre possiblement fractionnaire. 
Rappelons que la fonction $L_S(s,\pi,\rho)$ ne dépend que du $L$-paquet qui contient $\pi$.
Elle est donc par définition stable.

La formule des traces exprime \eqref{1.5} comme une somme sur des classes de conjugaison (stables) 
et la difficulté sera l'analyse de l'égalité qu'elle donne. À nos fins il est mieux de prendre 
d'abord une représentation $\rho$ irréductible et non triviale du groupe ${}^LG=\hat G\rtimes\Gal(K/F)$,
l'extension galoisienne étant suffisamment grande. Par exemple, pour $G=\{1\}$, $\rho$
est n'importe quelle représentation complexe, irréductible et de dimension finie
du groupe du Galois. En ce moment,
il ne semble pas que la question des zéros des fonctions $L(s,\pi,\rho)$ sur la droite
$\Re s=1$ doit nous préoccuper. Ils sont censés ne pas exister. Nous avons déjà constaté que les 
contributions à \eqref{1.5} ou \eqref{1.6} des représentations $\pi$ pour lesquelles le paramètre 
$\phi$ n'est pas trivial puissent se traiter  par récurrence, la dimension du centralisateur de 
$\phi$ étant plus petite que la dimension de ${}^LG$. On obtiendra à la
fin une formule dans laquelle toutes les 
représentations automorphes qui interviennent sont de type Ramanujan, de façon qu'en principe les 
fonctions $L$ pour lesquelles les sommes \eqref{1.5} ou \eqref{1.6} n'ont pas de pôle à la droite de la ligne 
$\Re s=1$. On obtient l'ordre du pôle en multipliant l'une ou l'autre des deux expressions par 
$s-1$ et en faisant $s\searrow1$. Il y a aussi d'autres façons analytiques d'extraire l'ordre du 
pôle de (1.5) ou (1.6) qui puissent s'avérer plus utiles, mais ce qui importe ici, c'est que ces 
informations sont recelées dans ces deux expressions. 

Si $\pi$ est de type Ramanujan on s'attend, et pour de très bonnes raisons, à ce que l'ordre 
$\mu(\pi)=\mu(\pi^{\st})$ du pôle de la dérivée logarithmique de $L(s,\pi,\rho)$ en $s=1$ soit la 
multiplicité de la représentation triviale dans la restriction de $\rho$ à ${}^\lambda H_\pi$,  de 
sorte que la somme
\begin{equation}\label{1.7}
\sum_{\pi^{\st}}\mu(\pi)m(\pi^{\st})\prod_{v\in S}\tr\pi_v^{\st}(f_v)
\end{equation}
sur des représentations de type Ramanujan se retrouve dans la formule des traces et peut en 
principe en être extraite. Si $\rho$ est irréduct\-ible et non trivial, cette multiplicité sera 
nulle sauf pour quelques $\pi$ tels que ${}^\lambda H_\pi\not={}^LG$. Par conséquent, elle sera nulle 
pour la plupart des $\pi$.  Les autres $\pi$ proviennent des groupes ${}^\lambda H_\pi$, donc des 
groupes $H$ qui seront tous de dimension plus petite que celle de $G$ et donc en principe bien 
compris, mais seulement une fois la fonctorialité établie. 
En utilisant la formule des traces pour ceux-ci et la restriction $\rho_H$ de $\rho$ à 
${}^LH\rightarrow {}^LG$, nous pouvons, de nouveau en principe, calculer (1.7) comme une somme sur 
$H$. Si la multiplicité de la représentation triviale dans $\rho_H$ est $\mu(\rho_H)$, la somme 
(1.7) doit être égale à 
\begin{equation}\label{1.8}
\sum_H \mu(\rho_H)\sum_{\pi_H^{\st}}m(\pi_H^{\st})
\prod_{v\in S}\tr\pi_v^{\st}(f_v). 
\end{equation}
Cette somme est renfermée dans la formule des traces pour les groupes $H$ rattachés aux groupes 
${}^\lambda H$ quoiqu'il faudra tenir compte des ennuis qui accompagnent les cha\^ines 
d'inclusion ${}^LG\supset{}^\lambda H\supset\dots$ 

Ce que nous venons de décrire serait une confirmation de la fonctorialité en utilisant la formule 
des traces stable, mais ce n'est pas cela que nous avons proposé. Nous proposons plutôt de vérifier la 
fonctorialité par les mêmes arguments. Cela ne s'avère pas facile, car l'analyse des sommes sur des 
classes de conjugaison qui intervienent dans la formule des traces ne l'est point (\cite{L3}). Une 
autre possibilité est suggérée par le lemme de la prochaine section. Quoique les logarithmes ont 
l'avantage formel important que les multiplicités $\mu(\pi)$ y apparaissent 
linéairement et sans s'encombrer de facteurs inutiles, traiter les 
fonctions $L$ elles-même semble plus facile du point de vue
analytique. En plus des difficultés analytiques elles-mêmes, dans l'égalité de \eqref{1.7} et 
\eqref{1.8} est cachée une difficulté même très grave, à savoir la fonctorialité pour les plongements
${}^\lambda H\rightarrow {}^LG$. Ici intervient en plus une difficulté mineure, que nous discuterons
lorsque l'occasion se présentera. Le groupe ${}^\lambda H$ plongé dans ${}^LG$ n'est pas 
lui-même dual à un groupe $H$.
Il est l'image d'un groupe ${}^LH$ par rapport à une application admissible surjective. Cela n'a
aucune importance, mais il faut l'expliquer. La difficulté clé est que la stratégie que nous
proposons exige que nous démontrions, en principe en utilisant encore la formule
des traces, la fonctorialité pour l'application ${}^LH\rightarrow {}^LG$. Jusqu'à présent nous
ne savons pas comment le faire. Nous abordons le problème à pas comptés. Expliquons la stratégie.

Définissons l'opérateur de Hecke $K_v^{\rho,(n)}$ de sorte que
\begin{equation}\label{RL.1.2}
\tr(\pi_v(K_v^{\rho,(n)}))=\tr\rho^{(n)}(A(\pi_v)),
\end{equation}
$\rho^{(n)}$ étant le produit symétrique de $\rho$ de degré $n$. Posons
$$
\bbL_v(s,\rho)=\sum_{n=0}^\infty q_v^{-ns}K_v^{\rho,(n)}.
$$
Bien que les sommes et produits que nous allons introduire par la suite convergent pour $\Re s$ 
suffisamment grand, il y a certains avantages à les traiter comme des séries formelles en 
$t=q^{-s}$, où $\N\mathfrakp_v=q_v=q^{\deg v}$ ou,
pour le cas des corps de nombres, comme des séries de Dirichlet formelles. Évidemment
\begin{equation}\label{RL.1.3}
\tr(\pi_v(\bbL_v(s,\rho)))=L_v(s,\pi_v,\rho)
\end{equation}
si $\pi_v$ est non ramifié
et $0$ sinon.

Notre but dans cet article et ceux qui le suivront
est d'entamer une discussion des conséquences possibles d'une formule
des traces stable. Nous sommes prêts à escamoter quelques questions de
base, qui pour les problèmes avec lesquels nous
commençons ne posent pas de difficulté. Dans une théorie systématique on dit qu'une représentation
$\pi_v$ irréductible de $G(F_v)$ est non ramifiée si le groupe $G_v$
est quasi-déployé sur $F_v$ et déployé sur une extension
non ramifiée et s'il y a un vecteur non nul stabilisé par un sous-groupe
hyperspécial donné. On dit qu'un paquet $\pi^{\st}_v $ est non ramifié si chacun de ses éléments
est non ramifié pour un choix convenable de ce sous-groupe hyperspécial. À partir du groupe $G$ sur $F$,
ou plutôt à partir d'une famille d'équations qui le définit,
on obtient des groupes non seulement sur chaque $F_v$ mais
aussi sur les anneaux $\calO_v$, au moins pour presque toute place finie
$v$.
Les sous-groupes $G(\calO_v)$ sont hyperspéciaux presque partout. Nous
exigeons en choisissant $S$ que ceci est le cas en-dehors de $S$. Si l'on choisit une autre famille
d'équations alors les sous-groupes $G(\calO_v)$ ainsi obtenus sont les mêmes presque partout. 
Ces choix faits, 
l'équation \eqref{RL.1.2} sera valable en-dehors de $S$ mais 
pour un seul élément du paquet des $\pi_v$ qui y interviennent.
Cet élément et l'algèbre de Hecke sont fixés par le choix de $G(\calO_v)$. Pour les
autres éléments du paquet, $\tr(\pi_v(K_v^{(n)}))=0$. Donc les formes précises de \eqref{RL.1.2}
et de \eqref{RL.1.3} sont
$$
\tr(\pi^{\st}_v(K_v^{(n)}))=\tr\rho^{(n)}(A(\pi_v))
$$
et
$$
\tr(\pi^{\st}_v(\bbL_v(s,\rho)))=L_v(s,\pi_v,\rho).
$$

La factorisation
$$
L_S(s,\pi,\rho)\prod_{v\in S}\tr(\pi^{\st}_v(f_v))=
\prod_{v\notin S}\tr(\pi_v(\bbL_v(s,\rho)))\prod_{v\in S}\tr(\pi^{\st}_v(f_v)).
$$ 
implique que la somme
\begin{equation}\label{1.9}
\sum_{\pi^{\st}}m(\pi^{\st}) L_S(s,\pi,\rho)\prod_{v\in S}\tr(\pi^{\st}_v(f_v))
\end{equation}
est donnée par la formule des traces stable pour la fonction
\begin{equation}\label{f}
f=\otimes_{v\notin S}\bbL_v(s,\rho)\bigotimes \otimes_{v\in S} f_v.
\end{equation}
Il y a dans ces formules un mélange d'additif et de multiplicatif. Néan\-moins on peut espérer 
pouvoir écarter encore et de la même façon par récurrence les contributions des représentations 
qui ne sont pas de type Ramanujan. La somme qui restera sera égale à la différence entre une somme 
donnée par la formule des traces pour le groupe $G$ de départ et pour des groupes reliés à des 
plongements non triviaux de $\SL(2)$ dans ${}^LG$. Dans la cinquième section nous abordons le cas 
le plus simple, celui du plongement principal de $\SL(2)$ dans ${}^LG$ et nous montrons comment 
déceler dans la formule des traces stable pour $G$ les représentations rattachées à celui-ci. 
Dans cet article, nous 
n'irons pas plus loin et nous sommes plus convaincus par l'élégance et la 
simplicité des principes nouveaux que par nos résultats. Nous espérons revenir dans un très proche 
avenir à ces questions, donc de dégager dans la formule des traces les contributions des autres 
représentations qui ne sont pas de type Ramanujan, donc rattachées à des $\phi$ qui ne sont ni 
triviaux ni principaux.

Cela étant fait, on arrivera à partir de la formule des traces stable à une expression, mettons 
$\Xi(s)$, pour une somme comme \eqref{1.9} mais dans laquelle seules les représentations de type 
Ramanujan interviennent, au moins en 
principe. Pour celles-ci, les fonctions $L_S(s,\pi,\rho)$ ne devraient avoir ni pôle 
ni zéro dans le domaine $\Re s>1$ et l'on espère pouvoir tirer cette conclusion de cette expression. De toute 
façon, les représentations $\pi^{\st}$ qui contribuent à cette expression seront de types différents, 
à savoir celles qui proviennent des sous-groupes ${}^\lambda H$ tels que la restriction de $\rho$ à 
${}^\lambda H$ ne contient pas de représentation triviale et les autres. Pour un $\rho$
donné ces premières n'ont pas d'influence sur le pôle de $\Xi(s)$ en $s=1$
et nous pouvons les jeter au rancart. 

Les classes qui proviennent par transfert d'un groupe $H$ tel que la composition de $\rho$ 
avec l'homomorphisme ${}^\lambda H\rightarrow {}^L G$ contient la représentation triviale sont les
seules qui
contribuent à cette partie principale. Ces groupes $H$ sont de dimension plus petite que celle de 
$G$, et, par con\-sé\-quent, en principe compris. On peut donc, à partir de la formule 
des traces stable pour ceux-ci, 
calculer leur contribution à la partie principale de $\Xi(s)$ en $s=1$. En sommant ces parties 
principales sur tous les $H$ et en comparant cette somme à la partie principale de $\Xi(s)$, on 
devrait trouver une égalité qui confirmerait encore une fois la fonctorialité, ou mieux mènerait à 
sa démonstration en général. Nous abordons ces calculs dans la cinquième section 
d'une façon modeste mais le vrai 
travail reste à faire. 

Il y a une difficulté au coeur de cette stratégie que n'est pas
encore tout à fait résolue mais qui ne nous semble pas grave. Nous utilisons
les fonctions $L$ et non pas leurs dérivées logarithmiques. Donc si
$\pi_H$ n'est pas l'image par fonctorialité d'une représentation d'un groupe $H'$
de dimension plus petite que celle de $H$ et si
la restriction de $\rho$ à ${}^\lambda H$ est la somme de $m$ fois la représentation triviale
et d'une représentation $\tau$ qui ne contient pas la représentation triviale, alors
en principe
$$
L_S(s,\pi,\rho)=\zeta_F^m(s)L_S(s,\pi,\tau),
$$
où le deuxième facteur n'a ni zéro ni pôle en $s=1$. Il a néanmoins une influence
sur la partie principale en $s=1$. Par contre, la formule des traces stable est une somme (!) 
sur les classes stables des $\pi$ qui ne mêle pas multiplicativement les fonctions
$L$ rattachées aux classes $\pi^{\st}$ différentes. 

La grosse difficulté sera l'analyse du comportement asymptotique de la somme \eqref{1.9} mais avec tous les 
termes enlevés pour lesquels $\phi$ est non trivial. Nous ne l'avons pas encore entamée,
mais les sections 3,4 et 5 sont préparatoires à cette fin.

\section{Le cas des corps de fonctions}

L'opérateur $\prod_{v\notin S} \bbL_{v}(s,\rho)$ admet une interprétation fort agréable dans le cas
où le corps global $F$ est celui des fonctions rationnelles d'une courbe $X$ 
géométriquement connexe définie
sur un corps fini $\kappa=\FF_{q}$. L'ensemble fini $S$ est dans ce cas un sous-schéma fermé 
d'une courbe projective lisse $X$ sur $\kappa$. Soit $U$ son complément.

En développant le produit infini $\prod_{v\notin S} \bbL_{v}$ puis en regroupant les
termes ayant le facteur $q^{-ds}$ pour chaque entier naturel $d$, on a l'identité de séries formelles
\begin{equation}\label{serie}
\prod_{v\notin S} \bbL_{v}(s,\rho) = \sum_{d\in \NN} q^{-ds} 
\sum_{\sum_{i} d_{i}v_{i}\in U_{d}(\kappa)}
\prod_{i} K_{ v_{i}}^{\rho,(d_{i})}
\end{equation}
où $U_{d}$ désigne la $d$-ième puissance symétrique de $U$ et où la deuxième sommation est étendue
sur l'ensemble des diviseurs effectifs de degré $d$ de $U$. Les opérateurs
\begin{equation} \label{bbK}
\mathbb K_{\rho,d}=\sum_{\sum_{i} d_{i}v_{i}\in U_{d}(\kappa)}
\prod_{i} K_{ v_{i}}^{\rho, (d_{i})}
\end{equation}
admettent une interprétation géométrique.

Au lieu du groupe $G$ défini sur $F$, on se donne un $X$-schéma en groupes
lisse de fibres connexes ayant $G$ comme fibre générique et dont la restriction à l'ouvert $U$ est 
réductif, même quasi-déployé et déployé sur un recouvrement étale et fini. Il est commode d'utiliser
la même lettre $G$ pour désigner ce schéma en groupes réductifs.
Nous ne faisons en effet que répéter dans un langage géométrique
nos conditions sur l'ensemble $S$.  Soit $\Oc_{x}$
le complété formel de $\Oc_{X}$ en le point $x$.
En une place $x\notin U$, la donnée de ce schéma en groupes fixe un sous-groupe
compact $G(\Oc_{x})$ de $G(F_{v})$ qui n'est pas nécessairement maximal.
Le choix de ce schéma en groupes est essentiellement équivalent au choix 
d'un sous-groupe compact ouvert 
$G(\Oc_{\bbA})$ du groupe adélique $G(\bbA)$. L'espace des doubles classes
$$G(F)\backslash G(\bbA) /G(\Oc_{\bbA})$$
peut s'interpréter  en termes de $G$-torseurs au-dessus de $X$. Le champ ${\rm Bun}_{G}$
classifiant des $G$-torseurs sur $X$ est un champ algébrique d'Artin 
localement de type fini.\footnote{Pour plusieurs raisons
un fondement des faisceaux $\ell$-adiques sur le champ de modules
des G-torseurs nous manque. D'une part, il manquait une théorie adéquate des
faisceaux $\ell$-adiques sur un champ algébrique d'Artin. Cette théorie 
a été heureusement établie par Laszlo et Olsson (\cite{LO}). D'autre part, 
les spécialistes de la théorie dite de Langlands
géométrique se sont
intéressés à la construction des formes cuspidales particulières. À leurs fins
un fondement n'était pas strictement nécessaire.
Puisque nous nous intéressons à la formule des traces, donc à toutes
les formes automorphes
simultanément, un fondement adéquat nous est indispensable. Nous espérons le mettre
en place dans un article postérieur. Pour le moment deux références
possibles sont \cite{LMB} et le projet
{\tt{http://www.math.upenn.edu/$\sim$kresch/teaching/stacks.html}}.
Pour le champ ${\rm Bun}_{G}$ lui-même on peut consulter un article de
Jochen Heinloth \cite{H}.}
L'ensemble de ses $k$-points est une réunion
disjointe des espaces des doubles classes 
\begin{equation}\label{A.1}
	{\rm Bun}_{G}(\kappa)= \bigsqcup_{\xi\in {\rm ker}^{1}(F,G)}
	G_{\xi}(F)\backslash G_{\xi}(\bbA) /G_{\xi}(\Oc_{\bbA})
\end{equation}
pour une collection de formes $G_{\xi}$ de $G$, les indices
$\xi$ parcourant le sous-ensemble ${\rm ker}^{1}(F,G)$ des classes $\xi\in {\rm H}^{1}(F,G)$
localement triviales pour la topologie étale, 
la forme $G_{\xi}$ étant construite à partir de la torsion intérieure associée
à $\xi$.

Les opérateurs de Hecke sont incarnés par certains faisceaux pervers $\ell$-adiques
sur le champ des modifications.
On note $\Hck$ le champ de modules des quadruplets $(x,E,E',\phi)$ où $x$ est un point de $X$,
où $E$ et $E'$ sont des $G$-torseurs sur $X$ et où 
$\phi$ est un isomorphisme entre les restrictions de
$E$ et $E'$ à $X-\{x\}$. Autrement dit, $\phi$ est une modification de $E$ en le point $x$ dont le 
résultat est $E'$. Il est raisonnable de se restreindre, au moins au début, aux lieux de
modification $\{x\}\subset U$. La fibre de ${\Hck}$ au-dessus d'un couple $(x,E)$ est alors
une Grassmannienne affine. C'est un ind-schéma muni d'une action de $G(\Oc_{x})$. Pour toute 
représentation de dimension finie 
du groupe dual $\rho: \,^{L} G \to {\rm GL}(V)$, on dispose d'un complexe $\mathcal A_{\rho}$
sur ${\Hck}_{U}$ tel que la restriction de $\mathcal A_{\rho}$ à la fibre 
de ${\Hck}$ au-dessus d'un point $x\in U$ et d'un $G$-torseur $E$ fixé est le faisceau pervers
sur la Grassmannienne affine qui correspond à $\rho$ par l'équivalence de Satake géométrique
(\cite{MV}).
Le support de ce faisceau est donc de dimension finie quoique la fibre elle-même est de
dimension infinie.
Ces faisceaux pervers $\mathcal K_{\rho}$ s'organisent en un complexe au-dessus de $\Hck$.
Ce complexe est un faisceau pervers après un bon décalage. L'opérateur sur $D^{b}(\Bun_{G})$ 
\begin{equation}\label{K1}
\calF \mapsto \mathbb K_{\rho}(\calF)=\pr_{2,!} (\pr_1^* \calF \otimes \mathcal K_{\rho})
\end{equation}
est l'interprétation de $\bbK_{\rho,1}$ de la formule \eqref{bbK}
dans le cadre des faisceaux. Ici, $\pr_{1}$
et $\pr_{2}$ sont les projections de $(x,E,E',\phi)$ sur la composante $E$ et $E'$ respectivement.

Grâce à \eqref{A.1} et à la trace de Frobenius
sur les fibres 
$$
\sum_i(-1)^i\tr({\rm {Fr}},H^{(i)}(\calF_x))
$$
qui permet le passage, décrit par le dictionnaire de Grothendieck, d'un complexe de faisceaux à une fonction, 
la théorie
classique des formes automorphes se transforme en une théorie géométrique,
sauf que le corps des nombres complexes est remplacé par la clôture algébrique
de $\QQ_\ell$ (\cite{Lm2}).
Rappelons que l'opérateur de Hecke habituel admet l'interprétation géométrique
$$
\calF \mapsto \Hck_{\rho}(\calF)=\pr_{\Bun_{G} \times U,!} (\pr_1^* \calF \otimes \mathcal K_{\rho})
$$
où $\pr_{\Bun_{G} \times U}$ est la projection $(x,E,E',\phi) \mapsto (E',x)$.
Ainsi l'opérateur \eqref{K1} consiste à intégrer $\Hck_{\rho}(\calF)$ le long de $U$.

Pour que
le paramètre de Arthur 
\begin{equation}\label{RL.2.1}
\sigma=\phi\times\psi:SL_{2} \times  W_{F}\to \,^{L} G
\end{equation}
définisse un faisceau $\ell$-adique il faut interpréter le groupe ${}^LG$
comme un groupe défini sur la clôture algébrique $\bar\QQ_\ell$ du corps $\QQ_\ell$ et $W_F$ comme un sous-groupe
du groupe de Galois $\Gal(F^{\sep}/F)$, le sous-groupe des éléments dont l'image dans $\Gal(F^{\nr}/F)$
est une puissance de $\mathbb F$. On rattache d'abord à $\sigma$ et $\rho$
une gradation de l'espace de $\rho$, définie par l'action du tore diagonal de $\SL_{2}$, donc par les valeurs
propres $\{m_j,m_j-2,\dots,-m_j\}$
de $\sigma_j(1,-1)$, où $\sigma_j$ est donné
dans \eqref{RL.1.0} 
et ensuite un système local gradué, donc une représentation $\psi'$ du groupe $W_F$
compatible avec la gradation. Elle est
donnée
par
$$
\psi':\,w\to \phi\left(\begin{matrix} q^{n/2}&0\\0&q^{-n/2}\end{matrix}\right)\psi(w),
$$
où l'image de $w$ dans $\Gal(F^{nr}/F)$ est $\FF^n$. Nous soulignons que nous supposons
que $\overline\QQ_\ell$ est plongé dans $\CC$ pour que $q^{n/2}$ puisse s'interpréter
comme un élément de $\overline\QQ_\ell$. Il y a donc une ambiguité curieuse quoique
familière de signe. 
 
Le faisceau $\calF$ est un faisceau propre pour les opérateurs de Hecke ayant comme 
valeur propre le paramètre de Arthur \eqref{RL.2.1}
s'il existe un isomorphisme\footnote{Le lecteur plus à son aise avec la
théorie classique est encouragé à s'arrêter sur l'équation \eqref{K1}
et se convaincre qu'elle implique la définition classique de Hecke pour les fonctions
rattachées aux faisceaux et que la relation qui suit implique que la fonction
rattachée au faisceau par la trace satisfait, grâce à l'isomorphisme de \cite{MV}, 
à celle connue depuis Hecke,
sauf que les coefficients appartiennent à un corps différent.}
$$
\Hck_{\rho}(\calF)\simeq\calF\boxtimes \mathcal L_{\rho}(\sigma)
$$
où $\mathcal L_{\rho}(\sigma)$ est le système local gradué obtenu en composant
$\rho\circ \sigma: W_{F} \times SL_{2} \to\, ^{L} G$ puis en restreignant à $W_{F}$.
En particulier,
si la restriction de $\sigma$ à $\SL_{2}$ est triviale, $L_{\rho}(\sigma)$ est concentré
en degré zéro. De plus, on veut une compatibilité avec le produit tensoriel de représentations
de $^{L}G$.
Si $F$ est un faisceau propre pour les opérateurs de Hecke ayant comme 
valeur propre le paramètre de Arthur, on a 
$$
\mathbb K_{\rho}(F)\simeq F \otimes \rmH^{*}_{c}(U,\mathcal L_{\rho}(\sigma)).
$$

Pour tout $d\in\NN$, le champ de modification $\Hck_{d}$ au-dessus de la puissance 
symétrique $U_{d}$ classifie les quadruplets $(D,E,E',\phi)$ où $D\in U_{d}$ est un 
diviseur effectif de degré $d$ dans $U$, où $E$ et $E'$ sont des $G$-torseurs sur $X$ et 
où $\phi$ est un isomorphisme entre les restrictions de
$E$ et $E'$ à $X-D$. Au-dessus d'un $G$-torseur fixe $E\in \Bun_{G}$ et d'un diviseur 
effectif fixe $D=\sum_{i=1}^{r} d_{i} x_{i}$ de support réduit $x_{1}+\cdots+x_{r}$
de longueur $r$, $\Hck_d$ s'identifie au produit de $r$ copies de la Grassmannienne 
affine. Chacune correspond à un point $x_{i}$. Au-dessus de $\Hck_{d}$, on peut construire 
un complexe $\mathcal K_{\rho,d}$ tel qu'au-dessus d'un diviseur  $\sum d_{i}x_{i}$
se trouve $\bigotimes_{i=1}^{r}\mathcal K_{x_{i}}^{\rho,(d_{i})}$ où 
$\mathcal K_{x_{i}}^{\rho,(d_{i})}$ est le faisceau pervers
sur la Grass\-mann\-ienne affine correspondant à la $d_{i}$-ième puissance symétrique de $\rho$.
L'existence de ce faisceau pervers ne va pas de soi mais sera expliquée dans un prochain article
\cite{FN}.
Notons que dans le cas où $G=\GL_{n}$ et $\rho$ est la représentation standard, il a été construit
par Laumon dans \cite{Lm1} et la construction générale n'en est pas très différente.
L'opérateur sur $D^{b}(\Bun_{G})$ 
$$ 
F \mapsto  \mathbb K_{\rho,d}(F)= \pr_{2,!}(\pr_1^* F \otimes \mathcal K_{\rho,d})
$$
est l'interprétation de $\bbK_{\rho,d}$ de la formule \eqref{bbK}
dans le cadre des faisceaux. Ici, $\pr_{1}$
et $\pr_{2}$ sont les projections de $(D,E,E',\phi)$ sur la composante $E$ et $E'$ respectivement.

\begin{lemme}
Si $F$ est un faisceau propre pour les opérateurs de Hecke ayant comme 
valeur propre le paramètre de Arthur $\sigma$, on a 
\begin{equation}\label{RL.2.2}
\mathbb K_{\rho,d}(F)=F \otimes S^{d}\rmH^{*}_{c}(U,\mathcal L_{\rho}(\sigma)).
\end{equation}
\end{lemme}

\begin{proof}
Puisque ce lemme n'est donné qu'à titre d'arrière-fond de notre stratégie, 
nous nous contentons ici d'un argument fonctionnel et reportons la démonstration de l'égalité \eqref{RL.2.2}
de faisceaux
à l'article \cite{FN}. On note $f$ la fonction sur $G(F) \backslash G(\bbA)/G(\Oc_{\bbA})$
qui correspond au faisceau $F$ selon le dictionnaire de Grothendieck.
Par définition, $\mathbb K_{\rho,d}(f)$ est la somme sur les diviseurs effectifs $D$ de
degré $d$ dans $U$
$$ \sum_{D\in U_{d}(\mathbb F_{q})} K_{D,\rho}(f)$$
où l'opérateur $K_{D,\rho}$ attaché à un diviseur $D=\sum d_{i} x_{i}$ est le produit
$\sum_{i} K_{x_{i}}^{\rho,(d_{i})}$ dont le transformé de Satake est la fonction trace de
la représentation de $^{L} G$ sur la $d_{i}$-ième puissance symétrique $\rho^{(d_{i})}$. 

Puisque $F$ est un faisceau propre pour les opérateurs de Hecke avec valeur propre $\sigma$,
$f$ est une fonction propre pour chacun des opérateurs $K_{D,\rho}$ avec valeur propre
$$ \tr({\rm Fr}_{q}, \bigotimes_{i=1}^{s} S^{d_{i}} \mathcal L_{\rho}(\sigma)_{x_{i}}).$$
La formule des traces de Grothendieck-Lefschetz implique alors l'égalité
$$ \sum_{D\in U_{d}(\mathbb F_{q})}\tr(\sigma_{D}, 
\bigotimes_{i=1}^{s} S^{d_{i}} \mathcal L_{\rho}(\sigma)_{x_{i}})=\tr\left({\rm Fr}_{q}, 
\rmH^{*}_{c}(U_{d}\otimes_{\kappa} \bar \kappa,
S^{d} \mathcal L_{\rho}(\sigma)_{x_{i}})\right)
$$
qui est le reflet du lemme sur les fonctions.
\end{proof}

Dans le cas où $L_{\rho}(\sigma)$ est concentré en degré zéro, 
le complexe $\rmH^{*}_{c}(U\otimes_{\kappa}\bar \kappa,\mathcal L_{\rho}(\sigma))$
est concentré en degré $0,1$ et $2$ et en les seuls degrés $1$ et $2$ si $U$ est affine. 
Si $\mathcal L_{\rho}(\sigma)$ n'a ni sous-faisceaux ni quotients constants, alors
$\rmH^{*}_{c}(U\otimes_{\kappa}\bar \kappa,\mathcal L_{\rho}(\sigma))$
est concentré en degré $1$.
La puissance
symétrique $S^{d}\rmH^{*}_{c}(U,\mathcal L_{\rho}(\sigma))$ est calculée par la formule suivante
$$\bigoplus_{d_{0}+d_{1}+d_{2}=d} 
S^{d_{0}} \rmH^{0}_{c}
\otimes \bigwedge^{d_{1}} \rmH^{1}_{c}
\otimes S^{d_{2}} \rmH^{2}_{c} [-d_{1}-2 d_{2}]$$
où on a omis les parenthèses évidentes $(U\otimes_{\kappa}\bar \kappa,\mathcal L_{\rho}(\sigma))$ 
suivant $\rmH^{*}_{c}$. Au niveau des fonctions le décalage du degré n'est guère qu'une formalité.
Si $\mathcal L_{\rho}(\sigma)$ n'a pas ni sous-faisceaux ni de quotients constants, alors
$S^{d}\rmH^{*}_{c}(U,\mathcal L_{\rho}(\sigma))=0$ pour $d$ plus grand que $\dim 
\rmH^{1}_{c}(U\otimes_{\kappa}\bar \kappa,\mathcal L_{\rho}(\sigma))$. Cette dimension peut être calculée en
fonction de la ramification sauvage de $\mathcal L_{\rho}(\sigma)$ par la formule de 
Grothendieck-Ogg-Shafarevich \cite{R}.

\section{Hypothèses, mesures et la formule des traces stable}\label{hm}

\subsection{Hypothèses}\label{hypo}

Soit $G$ un groupe réductif sur un corps global $F$ qui est ou bien une extension finie de $\QQ$ 
ou bien le corps des fonctions rationnelles d'une courbe $X$ sur un corps fini $\kappa$ à $|\kappa|=q$ 
éléments. Notons $\bbA_F$ l'anneau des adèles de $F$.

La théorie des représentations automorphes se répartit en trois parties diffé\-rentes: 
l'endoscopie et les $L$-paquets qui ramènent tout à la formule des traces stable
sur des groupes quasi-déployés; la comparaison 
des représentations automorphes sur un groupe $G$ et celles sur sa forme intérieure 
quasi-déployée, donc la démonstration de théorèmes, dont celui dit parfois le théorème de 
Jacquet-Langlands est le plus simple, bien que lui appartient aussi à la théorie de
l'endoscopie; la fonctorialité pour les groupes quasi-déployés. 
Maintenant que le lemme fondamental est acquis, nous pouvons nous
attendre à des progrès rapides avec les problèmes des deux premières parties (\cite{A3}). Notre 
objectif principal dans cet article et ceux qui le suivront est d'entamer, si Dieu le permet, 
l'étude de la fonctorialité, qui sera de loin
la partie la plus difficile. Nous supposerons donc que $G$, et par conséquent aussi son 
groupe dérivé $G_{\der}$, est quasi-déployé.\footnote{Un rapporteur a été troublé par notre
tendance de louvoyer dans l'utilisation de cette hypothèse. La raison est simple. Du point
de vue d'une théorie ultime de la formule des traces stable, la formule pour un groupe $G$
arbitraire, quasi-déployé ou non, mais, disons, avec $G_{\der}$ simplement connexe, peut
être traitée comme une formule pour $G$ mais en utilisant comme un de ses groupes endoscopiques
sa forme quasi-déployé $G_{\qd}$ ou comme une théorie pour $G_{\qd}$ en utilisant
la transformée endoscopique $f^{G_{\qd}}$ de la fonction $f^G$ sur $G$. Notre hypothèse
est donc au fond inutile, mais elle rend parfois les explications plus simples.}

L'usage des $z$-extensions, notion qui a été formalisée dans \cite{K1}, permettra 
de se ramener au cas où le groupe dérivé est simplement con\-nexe. Une $z$-extension d'un 
groupe $G$ est, en particulier, une extension $G'\rightarrow G$ telle que le groupe dérivé 
$G'_{\der}$ est simplement connexe et telle que les homomorphismes $G'(F)\rightarrow G(F)$ et 
$G'(\mathbb A_F)\rightarrow G(\mathbb A)$ sont surjectifs. Dans un sens, la théorie des formes 
automorphes sur $G$ est contenue dans la théorie pour $G'$. Nous supposons donc dans cet article 
que le groupe dérivé $G_{\der}$ de $G$ est simplement connexe. Le cas d'un corps global de
caractéristique positive n'est pas traité dans \cite{K1}. Faute
de temps nous ne le traitons dans cet article
que d'une façon lacunaire. 

Soit $Z$ 
la composante
neutre du centre de $G$. Il y a une suite exacte
\begin{equation}\label{A}
\{1\}\rightarrow A\rightarrow Z\times G_{\der}\rightarrow G\rightarrow\{1\},\qquad
A=G_{\der}\cap Z.
\end{equation} 
Nous décrirons plus tard la structure des points sur $F$ dans ce que nous
appellerons la base de Steinberg-Hitchin mais qui ne signifie que l'ensemble des classes
de conjugaison stables. Le nom de Hitchin y est rattaché parce que nos réflexions
ont été influencées par la théorie pour les algèbres de Lie, mais les éléments
essentiels de la classification
des classes stables dont nous avons besoin sont antérieurs à cette théorie. Ils
se trouvent dans les articles de Steinberg et Kottwitz, quoique ni l'un ni l'autre
n'a insisté sur les paramètres linéaires. 
Il est possible de décrire la classification dans un cadre plus familier si  
le groupe fini $A$ est étale.  Si la caractéristique de $F$ est zéro,
c'est toujours le cas. Si l'ordre de $A$ est premier à la caractéristique
c'est aussi le cas, mais il n'est pas toujours ainsi.
Le cas le plus simple où cette difficulté ennuyeuse se présente est le groupe $G=GL(2)$
sur un corps de caractéristique $2$. 
Le groupe $A$ est le sous-groupe de $GL(1)$ défini par $z^2=1$. Pour simplifier la tâche
de rédaction nous supposons que $A$ est étale sur $F$. Cependant, pour créer une théorie complète 
il faudra enlever cette condition.  

Nous ne pouvons utiliser la formule des traces
effectivement que si le spectre automorphe contient une partie 
discrète. Ce n'est pas le cas si le centre de $G$ contient un $F$-tore déployé. Dans ce cas le 
quotient $G(F)\backslash G(\bbA_{F})$ a un volume infini. Il faudra alors ou bien remplacer 
$G(\mathbb A_F)$ par un groupe plus petit comme le fait Arthur ou bien remplacer $G(F)$ par un groupe 
plus grand. En fin de compte, la différence entre les deux choix est petite. Nous préférons 
toutefois remplacer $G(F)$ par un groupe plus grand, ce qui nous semble la solution la plus
élégante. 

Notons $Z_{\spl}$ le plus grand $F$-tore déployé dans $Z$  et $Z'_{\spl}$ le plus 
grand quotient de 
$G$ qui est un $F$-tore déployé. L'homomorphisme $Z_{\spl}\to Z'_{\spl}$ qui s'en déduit est alors 
une isogénie 
dont on notera $Z''_{\spl}$ le noyau. Notons $s_G$ le rang de $Z_{\spl}$ et de $Z'_{\spl}$. 
Choisissons un 
sous-groupe 
$$ \ZZ^{s_G} \subset Z_{\spl}(\bbA_{F})$$
tel que 
le quotient $\ZZ^{s_G} Z_{\spl}(F) \backslash Z_{\spl}(\bbA_{F})$ est compact.
L'intersection des deux groupes $\ZZ^{s_G}$ et $Z''_{\spl}$ est nécessairement triviale
de sorte que $\ZZ^{s_G}\rightarrow Z'_{\spl}(\bbA_F)$ est un plongement de $\ZZ^{s_G}$ dans 
$Z'_{\spl}(\bbA_F)$.
Pour tout sous-groupe $H$ de $G$, nous allons noter
$$H^{+}(F)=\ZZ^{s_G} H(F).$$
En particulier, on a $G^{+}(F)=\ZZ^{s_G} G(F)$. Le quotient 
$$\ZZ^{s_G} Z(F) \backslash Z(\bbA_{F})= Z^{+}(F) \backslash Z(\bbA_{F})$$
est un groupe compact. 

Le caractère central d'une représentation automorphe irréductible $\pi$ définit un caractère 
$\varpi$ de 
$Z(F) \backslash Z^{+}(F)=\ZZ^{s_G}$. Il existe un caractère $\varpi'$ de $Z'_{\spl}(F)\backslash 
Z'_{\spl}(\bbA_{F})$ dont la restriction à $\ZZ^{s_G}$ est $\varpi$. En remplaçant $\pi$ par 
$\pi\otimes {\varpi'}^{-1}$, on obtient une  représentation automorphe dont le caractère central 
est trivial sur $Z^{+}(F)$. On peut donc se restreindre aux représentations automorphes 
ayant cette propriété et ne considérer que le quotient
$$G^{+}(F) \backslash G(\bbA_{F})$$
qui a un volume fini pour la mesure invariante sur $G(\bbA_{F})$.

\subsection{Mesures}\label{mesure}

Il est
bien connu comment associer une mesure à une forme volume sur une variété différentiable 
réelle.
On peut faire de même pour une variété sur un corps local non-archimédien du moment qu'on 
a choisi une mesure sur celui-ci qui joue le rôle de la mesure de Lebesgue sur le corps des 
nombres réels. Il en est de même des adèles. Ces mesures peuvent se définir
d'une façon canonique qu'il est très important pour nous de comprendre mais qui est mal 
expliquée même dans 
les références les plus souvent citées (\cite{W1,W2}).
Puisqu'il nous sera vraiment important d'avoir une référence précise et aussi brève
que possible nous reprenons les définitions ici. On va donc commencer par fixer des mesures 
compatibles sur l'anneau des adèles $\bbA_{F}$ et sur les corps locaux $F_{v}$ pour toutes les 
places $v\in |F|$.

Une mesure invariante de Haar sur le groupe localement compact $\bbA_{F}$ est 
bien définie à une constante 
près. Comme le quotient de $\bbA_{F}$ par le groupe discret $F$ est un groupe compact, on peut 
normaliser la mesure d'une seule façon telle que le quotient $F\backslash \bbA_{F}$ ait la mesure 
un. C'est cette mesure invariante $dx$ sur $\bbA_{F}$ qu'on va choisir pour le reste de l'article,
mais il est préférable de ne pas la définir directement
par la condition que la mesure du quotient soit $1$ mais
à partir de mesures $dx_{v}$ sur les corps locaux $F_{v}$, $v$ une place 
de $F$, elles-mêmes définies à partir d'un caractère additif continu
global. C'est cette suite de définitions et sa logique qui est mal expliquée dans
\cite{W2} et qui mène à des difficultés sur le plan mnémonique et,
en fin de compte, sur le plan logique, lorsque l'on arrive à la formule des traces sur les groupes
réductifs. On commence avec un caractère global
$\chi$, qui donne à chaque place $v$ un caractère $\chi_v$. À partir de $\chi_v$
on définit une mesure $dx_v$. Si $|F|$ est l'ensemble de toutes ces places
alors $dx=\otimes_{v\in |F|} dx_{v}$ sera la mesure sur $\bbA_{F}$.
 
Considérons l'ensemble des caractères continus
$\chi:\bbA_{F} \to \CC^{\times}$ tels
que 
$\chi(b)=1$ pour tout $b\in F$ et $\chi(bb')=1$ pour 
tout $b\in F$ si et seulement si $b'\in F$. C'est un 
espace principal homogène sous le groupe $F^{\times}$. Il
est embêtant que l'existence d'un tel caractère est vérifiée pour les corps
de fonctions autrement que pour les corps de nombres et même peut-être troublant pour ceux qui sont épris de la
pierre de 
Rosette. Dans les deux cas n'importe quel caractère $\chi$ est défini
par ses composantes locales, $\chi(a)=\prod_v\chi_v(a_v)$

Si $F$ est le corps des fonctions rationnelles sur une courbe $X$ projective et
lisse sur le corps fini 
$\kappa$ à $q$ éléments, on choisit une forme différentielle méromorphe non nulle $\omega$ sur 
$X$. On a alors un caractère $\bbA_{F} \to \kappa$ défini par 
$$
x\mapsto \res(x \omega)=\sum_v\tr_{\FF_v/\kappa}(\res_v(x \omega)) 
$$
qui est 
trivial si $x=b\in F$. On en déduit un caractère $\chi:\bbA_{F}\to \CC^{\times}$ défini par
\begin{equation}\label{A.3.1}
x\mapsto \exp\left( \frac{2i\pi}{p}  \tr_{\kappa/{\mathbb F_{p}}}(\res(x\omega))\right).
\end{equation}
Si $F=\mathbb Q$ et $x\in\mathbb A_F$, on pose
$$
\chi_0(x)=\prod_v\chi_v(x_v),
$$
où $\chi_\infty(x_\infty)=\exp(-2\pi ix_\infty)$ et
$\chi_p(x_p)=\exp(2\pi ix')$ si $x'\in\mathbb Q$ est un nombre rationnel dont le dénominateur est une 
puissance de $p$ et tel que $|x'-x_p|_p\leq 1$. Enfin, pour n'importe quelle extension finie $F$ 
de $\mathbb Q$, on pose
\begin{equation}\label{A.3.2}
\chi(x)=\chi_0(\tr_{F/\mathbb Q}(x)).
\end{equation}

Ayant fixé un caractère $\chi$ de $\mathbb A_{F}$, on a pour toute place $v$ un caractère 
$\chi=\chi_v$ de $F_v$. Il existe alors une unique mesure de Haar $dx=dx_v$ sur $F_{v}$ qui est 
autoduale par rapport à la transformation de Fourier
$$
\hat f(y)=\int f(x)\chi(xy)dx,
$$
c'est-à-dire une mesure par rapport à laquelle on a $\hat{\hat f}(x)=f(-x)$. Si $v$ est non 
archimédien et 
$$
\mathfrak d_v^{-1}=\mathfrak d_v^{-1}(\chi)=\{x\,|\,\chi(xy)=1\,\forall\,y\in \mathcal O_v\},
$$
alors la mesure autoduale assigne la mesure $\N\mathfrak d_v^{-1/2}$ à $\mathcal O_v$. 
L'idéal $\mathfrak d_v$ est lié à la différente locale, mais il n'est pas égal à cette différent;
il est rattaché au caractère. À cause de lui il n'y a pas des calculs
locaux canoniques. L'idéal $\mathfrak d_v$ est toutefois égal à $\mathcal O_v$
presque partout. Donc il y a des formules
canoniques presque partout. Des calculs habituels montrent qu'à une place réelle la 
mesure autoduale rattachée au caractère $\chi(x)=\exp(2\pi iyx)$ est $|y|^{1/2}dx$; à une place 
complexe la mesure rattachée à $\chi(z)=\exp(2\pi i\Re(\bar wz))$ est $|w|dxdy$, où $z=x+iy$,
$w=u+iv$, 
$|w|=\sqrt{u^2+v^2}$. 

Ayant fixé le caractère global $\chi$, nous avons fixé en même temps les 
mesures autoduales locales $dx_v$ aussi bien que la mesure globale $dx=\prod_vdx_v$ sur $\mathbb A_F$
qui sera aussi autoduale.  
Cette dernière mesure est indépendante du choix de $\chi$ en vertu de la formule du produit 
$\prod_v|b|_v=1$ pour tout $b\in F^\times$. Puisque le groupe dual à 
$F\backslash\bbA_F$ est $F$ il suit facilement de la théorie générale
de la transformée de Fourier pour les groupes compacts que 
\begin{equation}
{\rm mes}(F \backslash \bbA_{F})=1
\end{equation}
et
que
$$
\sum_{b\in F} \hat f(-b)=\sum_{b\in F}f(b).
$$
C'est là la formule de Poisson qui donne l'équation fonctionnelle des
fonctions $L$ rattachées à $F$.

Soit $X$ une variété algébrique lisse de dimension $n$ sur $F$. La donnée d'une $n$-forme 
différentielle partout non nulle $\omega$ sur $X$ définit une mesure $|\omega|$ sur l'espace 
topologique $X(\bbA_{F})$. Pour tout point $x\in X(F_{v})$, il existe un voisinage analytique de 
$x$ dans $X(F_{v})$ isomorphe à un polydisque ouvert de coordonnées $a_{1},\ldots,a_{n}$. Sur ce 
polydisque la $n$-forme $\omega$ s'écrit $\omega=f da_{1}\wedge\cdots\wedge da_{n}$ où $f$ est une 
fonction analytique inversible sur le polydisque. La mesure $|f| da_{1} \cdots da_{n}$, transferée 
au voisinage de $x$, est en fait indépendante du choix des coordonnées locales. Cette mesure ne 
dépendant que de la $n$-forme $\omega$, on la note $|\omega|_{v}$. Si on se donne en plus un 
modèle lisse $\mathcal X$ de $X$ sur $\Oc_{v}$ et si on suppose que la $n$-forme $\omega$ s'étend 
en une $n$-forme invariante partout non nulle sur le schéma $\mathcal X$, on a la formule
\begin{equation}\label{Weil}
\int_{X(F_{v})} 1_{\mathcal X(\Oc_{v})}|\omega|_v= 
q_{v}^{-n}\N\mathfrak d_v^{-n/2} | \mathcal X(\kappa_{v})|
\end{equation}
où $\kappa_{v}$ est le corps résiduel de $\Oc_{v}$, $q_{v}$ est son cardinal, $| \mathcal 
X(\kappa_{v})|$ désigne le nombre de $\kappa_{v}$-points de $\mathcal X$ et finalement 
$1_{\mathcal X(\Oc_v)}$ désigne la fonction caractéristique du compact $\mathcal X(\Oc_v)$.
Le cas le plus simple de ces principes est évidemment le cas d'une forme invariante
sur un espace vectoriel de dimension finie.

Soit $G$ un groupe algébrique.\footnote{Bien que ces définitions
préliminaires sont valables pour tout groupe, nous supposons
dès le début qu'il est réductif. Il ne faut pas qu'il soit quasi-déployé
mais arrivés à la formule des traces stable, nous
posons cette condition supplémentaire. Rappelons que la théorie
de la formule stable, donc de l'endoscopie, ramène l'étude des formes automorphes sur un groupe
réductif arbitraire à la théorie pour les groupes quasi-déployés.
Cependant, nous ne supposons pas
qu'ils se déploient sur une extension non ramifiée de $F$. 
\vskip .1pc
\noindent
Il est très répandu
parmi les géomètres d'éviter la ramification ou de ne traiter
que la ramification modérée. D'un point de vue géométrique
il y a certainement des avantages en ne considérant que la théorie
des opérateurs de Hecke et en évitant maints problèmes de l'analyse
harmonique non invariante, ou de la ramification qui est
si répandue dans l'arithmétique. Ces problèmes éclairent toutefois
maintes questions de structure qui nous seront
importantes lorsque nous examinerons la formule des traces.}
La donnée d'une $n$-forme invariante 
sur $G$ est équivalente à la donnée d'un vecteur $\omega$ non nul dans le $F$-espace vectoriel
$\wedge^{n} \mathfrakg$ de 
dimension un où $\mathfrakg$ est l'algèbre de Lie de $G$. 
La $n$-forme invariante $\omega$
associée est alors non nulle partout sur $F_v$
et presque partout sur $\mathcal O_v$. Pour toute place $v$, elle définit une mesure invariante 
$|\omega|_{v}$ sur $G(F_{v})$. Il n'est en général pas possible de prendre la mesure 
$\bigotimes_{v\in |F|}|\omega|_{v}$ sur le groupe adélique car la formule \eqref{Weil}
implique que le produit infini 
$\prod_{v}\mu_{v}$ est en général divergent, où $\mu_{v}$ désigne la mesure  du sous-groupe compact 
$G(\Oc_{v})$:
$$
\mu_{v}=\int_{G(F_{v})} 1_{G(\Oc_{v})} |\omega|_{v}=\int_{G(\Oc_v)}|\omega_v|.
$$
En particulier, si $G$ est le groupe multiplicatif défini sur le corps des nombres rationnels $\mathbb 
Q$, pour tout nombre premier $p$, la mesure locale $\mu_{p}=1-p^{-1}$ est la valeur
en $1$ de l'inverse du facteur en 
$p$ de la fonction zêta de Riemann. 

La nécessité de modifier les mesures locales avant de prendre le produit pour arriver à une mesure
bien définie exige l'utilisation de plusieurs mesures locales. Bien qu'elle ne soit
pas absolument nécessaire, quelques auteurs introduisent non pas seulement la mesure produit
globale mais une renormalisation de cette mesure. Dans cet article nous utilisons
une notation qui distingue toutes ces mesures, en ajoutant des notations plus simples pour celles, qui
au nombre de deux ou trois, seront utilisées dans des articles
à suivre. Posons d'abord $d_{\geom}g_v=|\omega|_v$. Elle est la mesure donnée directement
par la forme $\omega$. Elle dépend évidemment de cette forme.
Soit $\sigma_{G}$ la représentation du module galoisien des caractères $G\to \Gm$ définis 
sur la clôture algébrique 
$\bar F$. Ils sont aussi définis sur la clôture 
séparable. La fonction $L$ d'Artin rattachée à la représentation sur $\sigma_G$
admet un dévelop\-pement en produit eulérien
$$L(s,\sigma_{G})=\prod_{v} L_{v}(s,\sigma_{G})$$
sur le demi-plan $\Re s>1$. La mesure locale normalisée sur $G(F_{v})$ est définie par 
\begin{equation}\label{mesure locale normalisee}
dg_v=d_{\norm}g_{v}= L_{v}(1,\sigma_{G})d_{\geom}g_{v}.
\end{equation}
Une fois les définitions
de base bien comprises, nous n'utiliserons pour cette dernière 
que la notation $dg_v$ sans indice inférieur supplémentaire
parce qu'elle est la mesure locale principale.

La mesure produit 
\begin{equation}\label{mesure produit}
dg=d_{\rm prod} g=\bigotimes_{v\in |F|} dg_{v}
\end{equation}
ne dépend pas  de $\omega$ à cause de la formule $\prod_{v\in |F|} |c|_v=1$ valable pour tout 
$c\in F^{\times}$. Cette mesure produit est bien définie d'après Ono (voir l'appendice 2 de 
\cite{O3}). Si $f$ est une fonction de la forme $f=\otimes_{v\in |F|} f_{v}$ avec 
$f_{v}$ lisse à support compact et égale à la fonction caractéristique
de ${G(\Oc_{v})}$ presque partout, le produit infini 
$$\prod_{v\in |F|} \int_{G(F_{v})} f_{v} dg_{v}$$
est absolument convergent.

On ajoute souvent une renormalisation supplémentaire globale,
mais nous préférons ne pas le faire. Néanmoins nous l'expliquons. Le rang de la partie triviale 
du module galoisien $\sigma_{G}$ est égal au rang $s_G$ de $Z_{\spl}$, et la fonction 
$L(s,\sigma_{G})$ a un pôle d'ordre $s_{G}$ en 
$s=1$. Si 
\begin{equation}\label{rho G}
\rho_{G}=\lim_{s\searrow 1} (s-1)^{s_G} L(s,\sigma_{G}),
\end{equation}
alors
\begin{equation}\label{mesure globale normalisee}
d\tilde g= \rho_{G}^{-1} d_{\rm prod} g
\end{equation}
est la mesure doublement renormalisée, mais nous ne l'utiliserons pas souvent.
Nous introduisons pourtant les deux mesures, la mesure produit, ``$\rm meas_{\rm prod}$'', définie
par $d_{\rm prod} g$, les mesures locales étant renormalisées, et
simplement ``$\rm meas$'' pour la mesure doublement renormalisée, non pas
parce qu'elle est plus importante mais parce qu'elle est plus répandue.
Cela pourrait entraîner dans cet article l'apparence fréquente du premier symbole 
qui est moins commode, mais en fait nous l'éviterons. L'avantage principal de la mesure produit est qu'elle
est plus facile à manier dans le cadre adélique, traditionel pour l'analyse harmonique sur $G$
et que nous introduisons aussi sur la base de Steinberg-Hitchin.  

Nous avons choisi le sous-groupe discret $G^+(F)$ en sorte que
le quotient $G^+(F)\backslash G(\mathbb A_F)$ est de mesure finie.
Sa mesure 
$$
\meas(G^+(F)\backslash G(\mathbb A_F)),
$$ 
qui apparaît dans
la formule des traces,
se calcule facilement \`a partir de la mesure
de Tamagawa $\tau(G)$
telle que d\'efinie et calcul\'ee dans l'article de Ono (\cite{O3}).
Introduisons la suite exacte
$$
\{1\}\to G_1\to G\to G_2\to\{1\}
$$
où $G_1$ est connexe et $G_2$ est un tore déployé de rang $s_G$.
L'image réciproque de $G_2(\mathbb A_F)$ 
dans $G(\mathbb A_F)$
est le groupe $G_A^1$ de \cite{O3}. 
La mesure
$$
\meas\left(G^+(F)\backslash G(\mathbb A_F)\right),
$$ 
calcul\'ee
par rapport \`a la mesure normalis\'ee globale, est le produit des mesures de
$G(F)\backslash G^1_A$ et du quotient de $G_A^1\backslash G_A$ par l'image
de $G^+(F)$, un groupe discret. Ce quotient $Q$ est compact, isomorphe \`a 
$\mathbb Z_{s_G}\backslash \mathbb R_{s_G}$
dans le cas d'un corps de nombres et \`a un groupe fini dans le cas d'un corps de fonctions.
Ces mesures sont 
calcul\'ees selon les d\'efinitions de \cite{O3}. Donc
\begin{equation}\label{A.3.4}
\meas\left(G^+(F)\backslash G(\mathbb A_F)\right)=\tau(G)\int_Qdt.
\end{equation}
Ce dernier facteur anodin $m_G=\int_Qdt$ d\'epend du choix
de $G^+(F)$, donc de $\ZZ^{\rm s_G}$. Le facteur $\tau(G)$ est le
nombre de Tamagawa.\footnote{Nous laisserons au lecteur le soin
d'expliquer pour un corps de
fonctions la signification g\'eo\-m\'et\-rique de $\mathbb Z^{s_G}$
aussi bien que la signification de $m_G$.}
\vskip 1pc
Si le groupe $G$ est un tore $T$, alors grâce à \eqref{Weil} les facteurs du produit
infini
\equation\label{prod}
\prod_{v}\int_{T(F_{v})} f_{v} |\omega_{v}|
\endequation
sont égaux presque partout à
$1/L_v(1,\sigma_T)$.
Si $T$ est anisotrope, c'est-à-dire si $s_T=0$,
alors le théorème de Dirichlet et ses généralisations impliquent que
le produit 
\equation\label{dir}
\prod_v\frac{1}{L_v(1,\sigma_T)}
\endequation
converge conditionellement si on prend les $v$ dans l'ordre qui correspond
à la taille de $q_v$, $q_v=|\FF_v|$, et \eqref{dir} est alors égal à $1/\rho_G$. 
Donc \eqref{prod} donne un résultat qui ne diffère de $\int fd\tilde g$
que du produit d'un nombre fini de facteurs locaux.

\subsection{La base de la fibration de Steinberg-Hitchin}\label{quotient adjoint}
Dans cet article nous ne traiterons que la partie elliptique
réguli\`ere de la formule des traces, qui sera utilisée pour une fonction
$f=\prod_vf_v$, les $f_v$ étant lisses et à support compact. Pour l'utiliser il nous
faudra une description convenable des classes de conjugaison 
semi-simples stables, donc de ce que nous appelons
la base de la fibration de Steinberg-Hitchin. Nous commencerons
avec le cas d'un groupe semi-simple simplement connexe et d\'eploy\'e,
pour passer ensuite au cas d'un groupe semi-simple et simplement
connexe quasi-d\'eploy\'e, lui-m\^eme suivi par le cas des groupes
semi-simples et simplement
connexes arbitraires, et enfin
au cas le plus g\'en\'eral que nous traitons, celui d'une $z$-extension, pour laquelle, nous le rappelons,
le groupe d\'eriv\'e est
simplement connexe. Nous soulignons toutefois
que stabilisation signifie implicitement un passage à un groupe quasi-déployé!

La partie principale de la formule des traces est une somme
sur les classes de conjugaison elliptiques et régulières;
la partie principale de la formule
des traces stable étant une somme sur les classes stables, elliptiques et régulières.
Leurs descriptions pour un tore et pour un groupe semi-simple, simplement
connexe sont différentes, mais, dans le cas général, il faut mêler les deux, car $G$
est construit à partir de $G_{\der}$ et du centre $Z$. À part quelques complications
cohomologiques qui restent à décrire les classes de conjugaisons semi-simples
dans $G(F)$, mettons stables, sont  
des produits d'un élément de $Z(F)$ et d'une classe de $G_{\der}(F)$.
L'analyse harmonique, locale ou globale, sur $G$ revient en fin de compte à l'analyse harmonique
sur ces deux facteurs. Il s'avère que les classes stables dans $G_{\der}(F)$
sont rattachées aux points d'un espace vectoriel de dimension finie sur $F$ et 
cela permet de traiter la somme qui apparaît dans
la formule des traces d'une tout autre façon que celles utilisées antérieurement. 

Pour les lecteurs avec une formation géométrique, nous anticipons nos conclusions
pour qu'ils soient bien conscients des problèmes analytiques que pose la formule des traces
et auxquels la structure additive que nous décrivons offre peut-être une solution,
et pour qu'ils ne soient pas trop entravés par leurs connaissances géométriques. Nous leur rappelons
surtout que le corps $F$ peut être aussi bien un corps de nombres algébriques qu'un corps
de fonctions!

La conclusion sera que la partie elliptique de la formule des traces
stable sera une somme sur $\eta\in\mathfrak h$,
où l'ensemble $\mathfrak h$ reste à décrire, de sommes sur l'ensemble 
\eqref{A.2.4}, un ensemble qui sera introduit
plus tard. Si nous divisons par $|A(F)|$, nous pouvons même faire
la somme sur $\mathfrak B_\eta(F)\times Z_\eta(F)$. Quoique la somme sur $\eta$ n'est pas
finie, pour une fonction $f=\prod_vf_v$ donnée, il n'y a qu'un nombre fini de $\eta$ pour
lesquels la contribution n'est pas $0$. Il suffit donc de traiter
la contribution d'un seul $\mathfrak B_\eta(F)\times Z_\eta(F)$. Il y a une action
simplement transitive du groupe $Z(F)$ sur $Z_\eta(F)$. La somme sur $Z_\eta(F)$
est donc une somme sur $Z(F)$ et par conséquent peut être traitée par la formule
de Poisson pour la paire $Z(F)\subset Z(\bbA_F)$ ou plutôt pour $Z^+(F)\subset Z(\bbA_F)$.
Ce qui est nouveau ici et qui n'est pas apparu auparavant, c'est --- à peu près --- que $\mathfrak B_\eta(F)$
est un espace vectoriel et que, pour un $z\in Z_\eta(F)$ donné, la somme sur $\mathfrak B_\eta(F)$
est une somme de Poisson des valeurs aux points de $\mathfrak B_\eta(F)$ d'une fonction adélique
dont le comportement est assez bon pour que la formule de Poisson puisse être utilisée. 
L'expression
``à peu près'' se rapporte aux conséquences des lemmes \eqref{getz} et \eqref{getz1}.
Pour utiliser la formule de Poisson un tronçonnage est nécessaire, car sinon on ne serait pas
en état de vérifier que la somme duale qui apparaît dans cette formule converge. Pour cela il
faut des majorations qui ne sont disponibles que localement et même alors
pas encore disponibles sauf dans quelques cas particuliers qui 
seront donnés dans \cite{L5}. Les problèmes ne se posent
qu'après le tronçonnage permis par le lemme de Getz, qui est introduit dans cet article
surtout pour pouvoir formuler la proposition \ref{prodinf}, mais qui a des objectifs plus 
ambitieux.

Après avoir expliqué toutes les définitions, nous donnerons l'exemple simple de $GL(2)$ pour
bien mettre en évidence la partie de la formule des traces que nous proposons
utiliser dans la formule de Poisson.   

Soit $T$ un tore de $G$ suppos\'e d\'eploy\'e et, pour le moment,
semi-simple et simplement connexe. Soient $\alpha_1,\dots,\alpha_r$
les racines simples de $T$ par rapport \`a un ordre choisi de fa\c con
arbitraire. Rattach\'es \`a ces racines simples sont les poids
fondamentaux $\mu_1,\dots,\mu_r$ d\'efinis par $\mu_i(\hat\alpha_j)=\delta_{i,j}$.
Soit $\rho_i=\rho_{\mu_i}$ la repr\'esentation de poids maximal $\mu_i$
et soit $b_i(t)=\tr\rho_i(t)$. Les $b_i$ sont alg\'ebriquement
ind\'ependants sur $F$ et  
le quotient de $T$ par le groupe de Weyl est l'espace affine 
$\Spec F[b_1,\dots,b_r]$. Ce quotient est la base de Steinberg-Hitchin $\mathfrak A$
dans lequel $b_1,\dots,b_r$ seront nos coordonn\'ees pr\'efér\'ees.
Par contre sur $T$ nos coordonn\'ees seront
$\gamma_i=t^{\lambda_i}$, o\`u
\begin{equation}\label{A.2.1}
\lambda_i=\sum_ja_{i,j}\mu_j,\qquad \mbox{avec}\qquad \det(a_{i,j})=\pm 1.
\end{equation}
Le choix pr\'ecis de la matrice d'entiers $(a_{i,j})$ est sans importance. Dans 
ce cas la base $\mathfrak A$
est aussi sa partie linéaire que nous notons $\mathfrak B$. Lorsque $G\neq G_{\der}$
cela n'est plus le cas. 

Le groupe fini des automorphismes du graphe de Dynkin se rel\`eve
\`a un groupe d'auto\-morph\-ismes, dont chacun est d\'efini sur $F$, du groupe semi-simple,
simplement connexe et d\'eploy\'e $G$.
Un groupe quasi-d\'e\-ploy\'e mais simplement connexe et semi-simple
se d\'efinit \`a partir d'un cocycle \`a valeurs dans ce groupe
d'automorphismes, donc d'un homomorphisme $\sigma\rightarrow\varphi(\sigma)$
du groupe $\Gal(\bar F^{\sep}/F)=\Gal_F$
dans le groupe d'automorphismes du graphe de Dynkin. Ce dernier groupe
agit aussi sur $T$, un sous-groupe de Cartan d\'eploy\'e, et sur $G$, 
les actions \'etant
d\'efinies sur $F$.  

Il y aura donc
une action \`a droite de $\Gal_F$ sur l'ensemble des re\-pr\'e\-sent\-ations
fondamentales $\rho_1,\dots,\rho_r$, 
$$
\sigma:\rho_i\mapsto \rho_{i^\sigma},\qquad
\rho_{i^\sigma}(g)=\rho_i(\varphi(\sigma)(g)),
$$
et par cons\'equent sur leurs
caract\`eres $b_1,\dots,b_r$. Observons que
$$
\rho_{(i^\tau)^\sigma}(g)=\rho_{i^\tau}(\varphi(\sigma)(g))=
\rho_i(\varphi(\tau)\varphi(\sigma)(g))
=\rho_i(\varphi(\tau\sigma)(g))=\rho_{i^{\tau\sigma}}(g).
$$
Soit provisoirement $G_\varphi$ le groupe quasi-d\'eploy\'e
d\'efini par le cocycle $\sigma\mapsto \varphi(\sigma)$.
L'action du groupe de Galois sur $G_\varphi(\bar F^{\sep})$,
qui comme ensemble n'est que $G(\bar F^{\sep})$, est 
$$
\sigma_\varphi:\,g\rightarrow \varphi(\sigma)(\sigma(g))=\sigma(\varphi(\sigma)(g)),
$$
$\sigma(g)$ \'etant d\'efini par rapport \`a $G$. 

Si $\{h_\sigma\}$ est, par rapport \`a cette action, n'importe
quel cocycle \`a valeurs
dans le centre de $G_\varphi(\bar F^{\sep})=G(\bar F^{\sep})$, alors 
$\rho_i(h_\sigma)\in \bar F^{\sep}$
est un scalaire $\epsilon_i(\sigma)$.
En plus, $\sigma(\epsilon_{i^\sigma}(\tau))\epsilon_i(\sigma)=\epsilon_i(\sigma\tau)$,
car
$$
\sigma(\rho_{i^\sigma}(h_\tau))\rho_i(h_\sigma)=
\rho_{i^\sigma}(\sigma(h_\tau))\rho_i(h_\sigma)=
\rho_i(\varphi(\sigma)(\sigma(h_\tau))h_\sigma)=
\rho_i(h_{\sigma\tau}).
$$
Par cons\'equent nous pouvons tordre l'espace lin\'eaire sur $F$
\`a coordonn\'ees $b_1,\dots,b_r$ en utilisant le cocycle
\begin{equation}\label{A.2.2}
\sigma: (b_1,\dots,b_r)\rightarrow
(b'_1,\dots,b'_r)\qquad\mbox{avec}\qquad b'_i=\epsilon_i(\sigma)b_{i^\sigma}.
\end{equation}
En particulier, en prenant pour $\{h_\sigma\}$ le cocycle trivial,
nous obtenons
une forme tordue de la base de Steinberg-Hitchin pour le groupe semi-simple,
simplement connexe et d\'eploy\'e qui est la base de Steinberg-Hitchin pour toute
forme int\'erieure du groupe quasi-d\'eploy\'e $G_\varphi$ rattach\'e \`a
l'homomorphisme $\varphi$, en particulier, pour $G_\varphi$ lui-m\^eme. Cette base
est toujours un espace vectoriel sur $F$.

Nous passons maintenant au cas g\'en\'eral où le groupe d\'eriv\'e $G_{\der}$
de $G$
est simplement connexe. 
Pour \'eviter des probl\`emes d'ins\'eparabilit\'e,
nous avons supposé 
que la caract\'eristique de $F$ est premi\`ere \`a l'ordre $|A|$
du groupe $A$ de \eqref{A}. Nous escamotons donc encore une fois une petite difficulté.

Selon la Proposition 2.2 de \cite{BT} tout \'el\'ement $g$ de $G(\bar F^{\sep})$,
donc en particulier de $G(F)$, est un produit $g=g_1g_2$,
$g_1\in G_{\der}(\bar F^{\sep})$, $g_2\in Z(\bar F^{\sep})$.
Si $g\in G(F)$, alors pour chaque $\sigma\in \Gal_F$,
on a $\sigma(g_1)^{-1}g_1\in A(\bar F^{\sep})$
et il est \'egal \`a $\sigma(g_2)g_2^{-1}$. Puisque $g_1$
n'est donn\'e qu'\`a partir d'un \'el\'ement de $A(\bar F^{\sep})$,
seulement l'image $\eta(g)$ du cocycle $\{\sigma(g_1)^{-1}g_1\}$
dans $H^{1}(F,A)$ est d\'efinie. Soit $\mathfrak h$ le 
sous-ensemble des $H^1(F,A)$ obtenu de cette
fa\c con.

Nous voulons d\'ecrire l'ensemble des points $F$-rationnels sur la base $\mathfrak A=\mathfrak A_G$ de
Steinberg-Hitchin de $G$ comme la r\'eunion sur $\eta\in\mathfrak h$, ou
m\^eme sur $H^1(F,A)$, des points rationnels
sur un ensemble qui est un quotient, à savoir,
\begin{equation}\label{A.2.4}
\{\mathfrak B_\eta(F)\times Z_\eta(F)\}/A(F)
\end{equation}
o\`u $\mathfrak B_\eta$ est un espace lin\'eaire, une forme
tordue de l'espace lin\'eaire sous-jacent \`a la base
de Steinberg-Hitchin de $G_{\der}$ sur laquelle $A(F)$ comme sous-groupe du centre
de $G_{\der}(F)$ agit et où $Z_\eta$ est le torseur sur $Z$
d\'efini par $\eta$ sur lequel $A(F)$ agit comme sous-groupe de $Z(F)$.
\`A cette fin, choisissons et fixons pour chaque $\eta\in\mathfrak h$ 
un cocycle $h=\{h_\sigma\}$ qui le repr\'esente. 
Si $\eta$ ou $h$ ne devient pas trivial dans $H^1(F,Z)$, alors $Z_\eta(F)$
est vide. Par contre, s'il devient trivial, alors $Z_\eta(F)$
peut \^etre identifi\'e avec l'ensemble des $g_2$ possibles pour ce
$h$ donn\'e.
L'ensemble de tout $z\in Z(\bar F^{\sep})$
tel que $\{\sigma(z)z^{-1}\}=\{h_\sigma\}$ sont les points d'un torseur
$Z_\eta$. C'est l'ensemble de tous les $g_2$ possibles pour
ce $h$. 

Supposons que 
\begin{equation}\label{A.2.5}
\sigma(g_1)^{-1}g_1=h_\sigma
\end{equation}
pour tout $\sigma$
et soient $b_1(g_1),\dots,b_r(g_1)$ les coordonn\'ees lin\'eaires
de l'image $\beta(g_1)$ de $g_1$ dans la base de Steinberg-Hitchin $\mathfrak A_{\der}$ 
de $G_{\der}$. 
Nous avons remarqu\'e que cette base est, pour n'importe quelle forme
tordue int\'erieure d'un groupe quasi-d\'eploy\'e semi-simple
et simplement connexe, \'egale comme vari\'et\'e sur $F$ \`a la base de 
Steinberg-Hitchin du groupe quasi-d\'eploy\'e dont nous avons d\'ecrit
la structure. En particulier, elle est un espace lin\'eaire, une forme
tordue de la base du groupe d\'eploy\'e. Nous tordons 
ce dernier espace lin\'eaire de dimension $r$ sur $F$ par le
cocycle rattach\'e \`a $\{h_\sigma\}$
comme dans l'\'equation \eqref{A.2.2}. V\'erifions que la condition
\eqref{A.2.5} est \'equivalente \`a la condition que $\beta(g_1)$
soit rationnel dans l'espace tordu. Nous soulignons que
l'espace lin\'eaire fondamental est alors la base de Steinberg-Hitchin
du groupe quasi-d\'eploy\'e!

Selon \eqref{A.2.2} la condition de rationalit\'e $b'_i=b_i$ pour tout $i$,
devient
$$
b'_i(\beta(g_1))=b'_i(g_1)=\epsilon_i(\sigma)\sigma(b_{i^\sigma}(g_1)).
$$
On a $\beta(\sigma(g_1))=\sigma(\beta(g_1))$. Plus
pr\'ecis\'ement, puisque la base de Stein\-berg-Hitchin de $G_{\der}$
est la forme tordue d\'ej\`a d\'ecrite de la base de Steinberg-Hitchin du groupe
semi-simple, simplement connexe d\'eploy\'e dont $G_{\der}$
lui-m\^eme est
une forme  tordue, on a
$$
b_i(\sigma(g_1))=\sigma(b_{i^\sigma}(g_1)).
$$
C'est-\`a-dire, l'action de $\Gal_F$ sur la base est
donn\'ee par \eqref{A.2.2}, le cocycle $\{\epsilon_i(\sigma)\}$
\'etant trivial, donc 
$$
\sigma:\,(b_i)\mapsto (\sigma(b_{i^\sigma})).
$$
Bien que l'action du groupe de Galois sur la base de Steinberg-Hitchin,
ou plut\^ot sur l'espace lin\'eaire y rattach\'e par les coordonn\'ees
$b_i$, soit celle sur la base de Steinberg-Hitchin du groupe d\'eploy\'e,
l'action sur $G_{\der}(\bar F^{\sep})$ est celle
d\'eduite de l'action sur $G(\bar F^{\sep})$, et non pas celle
d\'eduite de l'action sur un groupe d\'eploy\'e ou quasi-d\'eploy\'e.

Nous avouons que traîner toutes ces d\'efinitions
dans le bagage de l'article est fastidieux,
mais la structure lin\'eaire de la base de Steinberg-Hitchin
d'un groupe semi-simple et simplement connexe est \`a nos
fins un 
des \'el\'ements fondamentaux de la formule des traces. 

En traitant la classe des groupes r\'eductifs d\'ecrite 
dans la section \ref{hypo}, 
il faut
introduire le cocycle $\{\epsilon_i(\sigma)\}$
de \eqref{A.2.5} et la forme tordue de  la base de Steinberg-Hitchin du groupe
d\'eploy\'e
semi-simple, simplement connexe y rattach\'ee par \eqref{A.2.2}. 
Nous affirmons que par rapport \`a cette forme l'image $\beta(g_1)$
est d\'efinie sur $F$, donc 
que
$$
b_i(\beta(g_1))=b_i(g_1)=\epsilon_i(\sigma)\sigma(b_{i^\sigma}(g_1)).
$$
En effet,
$$
\sigma(b_{i^\sigma}(g_1))=b_i(\sigma(g_1))=b_i(h^{-1}_\sigma g_1)=
\epsilon^{-1}_i(\sigma)b_i(g_1).
$$

Nous ne pourrons pas \'eviter les ennuis qui proviennent du tore d\'eploy\'e
dans le centre de $G$. Alors les classes de conjugaison \`a traiter
seront non pas les
classes de $G(F)$ mais les classes de $G'(F)$. Puisque toute telle classe
contient le produit $zt$ d'un \'el\'ement de $\tilde Z$ et d'un \'el\'ement
de $G(F)$, la description de \eqref{A.2.4} reste {\it grosso modo} exacte. Il
faut simplement remplacer $Z_\eta$ par le groupe discret $Z^+_\eta=\ZZ^{s_G}Z_\eta$,
car si $z't'=g^{-1}zt g$, avec $t'$, $t$ et $g$ dans $G(F)$ et $z$, $z'$ dans 
$\ZZ^{s_G}$, alors $z^{-1}z'\in \ZZ^{s_G}\cap G(F)=\{1\}$.

Nous notons $\mathfrak A(F)$ l'ensemble des points sur $F$ sur la base de Stein\-berg-Hitchin de $G$,
avec une notation semblable pour $F_v$ ou pour $\bbA_F$.
Soit $\mathbf c(g)$ l'image dans $\mathfrak A$ d'un élément $g\in G$. 
Selon ce que nous venons de vérifier, il est une réunion sur $\eta$ de $\mathfrak B_\eta(F)\times Z^+_\eta(F)$ 
quotienté par $A(F)$. L'ensemble $Z_\eta$ est un torseur sur $F$ et donc un espace homogène principal
sur $Z^+(F)$ qui se plonge dans $Z(\bbA_F)$. D'autre part, $\mathfrak B_\eta(F)$ est un espace vectoriel
de dimension finie que se plonge dans $\mathfrak B_\eta(\bbA_F)$, un module libre de rang fini sur $\bbA_F$.
La technique qui vient à l'esprit pour exploiter la formule des traces
est alors d'utiliser la formule de Poisson pour $Z^+(F)\subset Z(\bbA_F)$ et pour
$\mathfrak B_\eta(F)\subset\mathfrak B_\eta(\bbA_F)$. Nous choisissons les mesures
globales et locales pour $Z$ et $\mathfrak B_\eta$
comme
dans la section \ref{mesure}.  

Toutes ces explications faites, il est évident que le cas de base est le cas d'un groupe semi-simple,
simplement connexe et dans les dernières sections de cet article nous nous bornerons à ce cas.
Il serait néanmoins utile d'ajouter quelques précisions pour le groupe $GL(2)$. La base de Steinberg-Hitchin
est de dimension deux et l'application $\mathbf c$ est donnée par $t\rightarrow (b,a)$ où
$X^2-bX+a$ est le polynôme caractéristique de $t$. Le paramètre $\eta$ est un choix $a_0$ de $a$ modulo 
l'ensemble des
carrés $c^2$ de $F^\times$, donc un choix du cocycle $\sigma\rightarrow \sigma(\sqrt{a_0})/\sqrt{a_0}$. 
On écrit alors $a=c^2a_0$, $b=d\sqrt{a_0}$,  de sorte que $d=b/\sqrt{a_0}$
est un point à coordonnées dans $F$ dans un espace principal homogène
sur le groupe additif. La dualité de Poisson additive sera utilisée pour la somme sur $b$ --- ou sur $d$.
Cette description
est valable même sur un corps de caractéristique $2$, mais alors la topologie convenable
sur $F$  n'est pas la topologie étale. Ce sont là des questions pour une autre occasion. 

\subsection{Les éléments elliptiques réguliers dans la formule stable des traces}
Pour cette partie de l'article, les travaux de Kottwitz et de Ono
inclus dans la bibliographie sont essentiels.\footnote{Puisque nous
voulons en principe traiter non seulement le cas où $F$ est un corps
de nombres mais aussi le cas où $F$ est un corps
de fonctions, ces références ne sont pas adéquates. 
Nous avons néanmoins décidé
de ne pas nous en occuper pour le moment.} Les notes
\cite{L2} pourraient aussi \^etre utiles.
Nous supposons toujours que le groupe $G_{\der}$ est quasi-d\'eploy\'e
et simplement connexe.

Rappelons que pour arriver \`a la partie elliptique r\'eguli\`ere
de la formule stable,
nous commençons avec la
somme
\begin{equation}\label{3.1}
\sum_{\gamma}\meas_{\rm prod}(T^+_\gamma(F)\backslash T_\gamma(\mathbb A_F))
\int_{T_\gamma(\mathbb A_F)\backslash G(\mathbb A_F)} f(g^{-1}\gamma g)d\bar g.
\end{equation}
La mesure du premier terme est donnée par $d_{\rm prod}t$. Celle
du second terme, $d\bar g$,
est la mesure quotient 
$$
\frac{d_{\rm prod}g}{d_{\rm prod}t}=\prod_vd\bar g_v\qquad \mbox{avec}
\qquad\bar g_v=
\frac{L(1,\sigma_G)}{L(1,\sigma_T)}\frac{d_{\rm geom}g_v}{d_{\rm geom}t_v}.
$$  
La vérification formelle de la formule \eqref{3.1} est facile. Le noyau
de notre op\'erateur est
$$
\sum_{\gamma\in T^+(F)}f(g^{-1}\gamma h).
$$
La trace s'obtient formellement de l'int\'egrale
$$
\int_{G^+(F)\backslash G(\mathbb A_F)}f(g^{-1}\gamma g)dg,
$$
qui est \'egale formellement \`a une somme sur les classes de
conjugaison $\{\gamma\}$ dans $G^+(F)$ de
\begin{equation}\label{3.2}
\int_{G^+_\gamma(F)\backslash G(\mathbb A_F)}f(g^{-1}\gamma g)dg.
\end{equation}
En fait cette somme n'a pas en g\'en\'eral de sens. N\'eanmoins
dans la vraie formule des traces (voir par exemple \cite{A2})
la somme des expressions sur l'ensemble des classes de conjugaison elliptiques et r\'eguli\`eres
dans $G^+(F)$ appara\^\i t. C'est cette somme qui nous int\'eresse dans
cet article, car toute recherche sur la formule des
traces commence avec l'\'etude de la contribution
des classes elliptiques r\'eguli\`eres. Si $\gamma$ est elliptique r\'egulier, 
alors son centralisateur
est un tore, donc connexe.\footnote{Pour le moment
nous ne donnons pas de r\'ef\'erence. Il faudra aussi trouver 
une d\'emonstration en caract\'eristique positive
du th\'eor\`eme 4.4 de Kottwitz (\cite{K1}).} 

Puisqu'il s'agit de la formule des traces stable nous travaillons
non pas avec des classes de conjugaison mais avec des classes
de conjugaison stable.
Selon \eqref{A.2.4} les paramètres de ces classes sont les points
dans la base de Steinberg-Hitchin $\mathfrak B_\eta(F)$.
Si $G$ n'est pas son propre groupe dérivé, la présence de $\mathbb Z^{m_{\rm sp}}$ 
et de $Z_{\eta}(F)$ rendent les formules plus compliquées. C'est pour cela que
nous avons supposé que $G=G_{\der}$, mais pour nous rappeler de temps en temps que ce n'est
pas le cas le plus général, nous utilisons parfois une notation qui tient compte du cas général. 
Ce qui est important c'est que $\mathfrak B_\eta$ est un espace vectoriel sur $F$.
 
Examinons les termes \eqref{3.2} en ajoutant
une somme sur $\ZZ^{s_G}$ pour que nous n'oublions pas que
ce sont les classes de conjugaison dans $G^+(F)$ qui comptent
et non pas celles dans $G(F)$, donc considérons
\begin{equation}\label{3.3}
\sum_{z\in\ZZ^{s_G}}\int_{G^+_\gamma(F)\backslash G(\mathbb A_F)}f(g^{-1}z\gamma g)dg.
\end{equation}
Le groupe $G_\gamma$, parfois d\'enot\'e $T_\gamma$, est le centralisateur
de l'\'el\'ement r\'egu\-lier semi-simple $\gamma$. 
Dans \cite{L2} la partie r\'eguli\`ere elliptique de la formule
des traces est convertie en une somme sur les groupes endoscopiques
de la formule des traces stables. Nous ne nous int\'eressons ici
qu'à la contribution de $G$ lui-m\^eme, qui est toujours, nous le rappelons,
tel que $G_{\der}$ soit simplement connexe, et en ce moment, pour simplifier
les choses, égal à $G_{\der}$. Nous expliquons bri\`evement
sa forme en citant \cite{K2} autant que possible,
car l\`a les notations et les explications sont plus
faciles \`a comprendre que celles de \cite{L2}.
\vskip1pc
\noindent
D'abord l'int\'egrale
$$
\int_{G^+_\gamma(F)\backslash G(\mathbb A_F)}f(g^{-1}z\gamma g)d_{\rm prod}g
$$
de \eqref{3.3}
est \'ecrite sous
la forme,
\begin{equation}\label{3.4}
\meas_{\rm prod}(T^+_\gamma(F)\backslash T_\gamma(\mathbb A_F))
\int_{T(\mathbb A_F)\backslash G(\mathbb A_F)} f(g^{-1}\gamma g)d_{\rm prod}\bar g.
\end{equation}
Ayant fait cette modification presque formelle, on prend d'abord
la somme de (3.3) non pas sur toutes les classes réguli\`eres $\{\gamma\}$ 
de $G(F)$
comme dans
(3.1) mais
sur tous les $\gamma$ qui sont conjugu\'es dans
$G(\mathbb A_F)$ \`a un $\gamma$ donn\'e. Cela donne un facteur,
d\'enot\'e dans \cite{K2} par 
$$
|\ker[\mathfrak E(T/F)\rightarrow\mathfrak E(T/\mathbb A_F)]|,\qquad\mbox{où} \qquad T=T_\gamma.
$$
Mais ensuite pour arriver \`a une somme sur des groupes
endoscopiques, on introduit un d\'eveloppement par rapport
aux caract\`eres d'un groupe $\mathfrak K(T/F)$. Cela m\`ene, au moins lorsqu'on 
ne consid\`ere que la contribution du groupe endoscopique
$G$ lui-m\^eme, donc la partie stable, au coefficient
\begin{equation}\label{3.5}
\iota(F,T,G)=\frac{|\ker[\mathfrak E(T/F)\rightarrow\mathfrak E(T/\mathfrak A_F)]|}
{\mathfrak K(T/F)},
\end{equation}
de \cite{K2}, la notation \'etant tout \`a fait pareille \`a celle de \cite{L2}.
Il est v\'erifi\'e dans \cite{K2}, Lemma 8.3.2, rédig\'e avant les d\'emonstrations
compl\`etes de la conjecture de Weil (\cite {K3})
sur les nombres de Tamagawa et de celle de Kneser sur le principe de Hasse, que
\begin{equation}\label{3.6}
\iota(F,T,G)=\frac{\tau(G)}{\tau(T)}.
\end{equation}

Introduire le facteur \eqref{3.5} et passer \`a la formule des traces stables
a pour cons\'equence que la somme sur $\gamma$ de (3.1) est remplac\'ee
par la somme sur les classes stables elliptiques r\'eguli\`eres, donc
effectivement par une somme sur les \'el\'ements elliptiques réguliers
de la base de Steinberg-Hitchin. Pour les sous-groupes de Cartan elliptiques
$m_T=m_G$ et, gr\^ace \`a \eqref{A.3.4}
le facteur de \eqref{3.4}  
$$
\meas_{\rm prod}(T^+_\gamma(F)\backslash T_\gamma(\mathbb A_F))=
m_T\rho_T\tau(T)=m_G\rho_T\tau(T)
$$
multipli\'e par \eqref{3.6} est \'egal \`a 
\begin{equation}\label{3.7}
m_G\tau(G)\rho_T=m_G\rho_G\tau(G)\frac{\rho_T}{\rho_G}.
\end{equation}
Ce qui reste est
$$
\int_{T(\mathbb A_F)\backslash G(\mathbb A_F)} f(g^{-1}\gamma g)d_{\rm prod}\bar g=
\prod_v\int_{T(F_v)\backslash G(F_v)}f_v(g^{-1}\gamma g)d\bar g_v.
$$
Nous introduisons la notation
$$
\int_{T(F_v)\backslash G(F_v)}f_v(g^{-1}\gamma g)d\bar g_v=
{\rm Orb}(\gamma,f_v),
$$
Prenant la somme de ces dernières intégrales sur des représentants des classes de conjugaison
dans la classe de conjugaison stable de $\gamma$, nous obtenons les
intégrales orbitales stables ${\rm Orb}(\gamma_{\st},f_v)$, $\gamma_{\st}$
dénotant la classe stable de $\gamma$. Dans la formule des traces ce sont les produits
\begin{equation}\label{3.8}
\prod_v{\rm Orb}(\gamma_{\st},f_v)
\end{equation}
qui interviennent. Avec une notation précise mais balourde on mettrait l'indice ``$\st$''
deux fois en l'ajoutant à ``$\rm Orb$''.

Le facteur $m_G\rho_G\tau(G)$ ne d\'epend pas
de $T$. Donc pour le moment nous le mettons au rancart. Dans lui sont cachés
des résultats cohomologiques et de mesures dont nous n'aurons plus de
besoin explicite. La repr\'esentation
$\sigma_G$ est une sous-repr\'esentation de $\sigma_T$. Soit $\sigma_{T/G}$
le quotient. Pour $T$ elliptique ce quotient ne contient pas la
repr\'esentation triviale et
$$
\frac{\rho_T}{\rho_G}=L(1,\sigma_{T/G})=\lim\limits_{s\searrow 1}L(s,\sigma_{T/G}).
$$
C'est le produit
\begin{equation}\label{3.9}
L(s,\sigma_{T/G})\prod_v{\rm Orb}(\gamma_{\st},f_v),\qquad s>1,
\end{equation}
qui nous int\'eressera dans la partie suivante. Ce produit est d\'efini
pour toute classe de conjugaison stable dans $T(F)$ m\^eme si elle n'est pas
elliptique. Toutefois, lorsque $s\mapsto1$, la limite n'existe pas
si la classe n'est pas elliptique. Puisque nous n'utiliserons jamais une seule
formule des traces mais toujours des sommes de plusieurs avec des signes différents,
il est tout à fait légitime d'ajouter une partie qui devient infinie lorsqu'on passe à 
la limite dans une des formules, pourvu que cette partie apparaisse avec un signe
opposé dans une autre. Nous aurons des exemples non pas dans cet article
mais dans la suite de celui-ci.

\subsection{Intégration le long des classes de conjugaison stable}

Il existe un ouvert dense $\mathfrak A^{\rs}$ de $\mathfrak A$ tel que pour tout $g\in G$ 
d'image $a\in \mathfrak A^{\rs}$, le centralisateur de $g$ est un tore. Pour tout élément
$a\in \mathfrak A^{\rs}(F_{v})$, l'ensemble des $F_{v}$-points de $G$ au-dessus de $a$,
non vide d'après Kottwitz \cite{K1}, forme une classe de conjugaison stable.  
L'ensemble $\mathfrak A^{\rs}(F_{v})$ est l'ensemble des classes de 
conjugaison stable semi-simples et fortement régulières.

Notre objectif dans cette section est  d'exprimer l'intégrale
$$ \int_{G(F_{v})} f_v(g_v) dg_v $$
en intégrant d'abord le long des fibres de $\mathbf c:G \to \mathfrak A$, puis en intégrant sur 
$\mathfrak A(F_{v})$. Nous aurons besoin de la formule qui en résulte
dans la dernière partie de l'article. Il s'agit d'abord de fixer des mesures de façon compatible.
Choisissons des mesures en suivant les
définitions de la section \ref{mesure}. En particulier, on fixe une forme volume 
$\omega_{G}$ sur $G$ et une forme volume $\omega_{\Af}$ sur $\Af$ définies sur $F$. 
Avec le choix du caractère additif, on en déduit des mesures $|\omega_{G}|_{v}$ sur $G(F_{v})$
et $|\omega_{\Af}|_{v}$. On a aussi des mesures normalisées $dg_{v}=L_{v}(1,\sigma_{G}) 
|\omega_{G}|_{v}$ et $db_{v}=L_{v}(1,\sigma_{Z}) |\omega_{\Af}|_{v}$, les facteurs de 
normalisation
$L_{v}(1,\sigma_{G})=L_{v}(1,\sigma_{Z})$ étant les mêmes.

La première observation est qu'il est possible de se restreindre 
à $\mathfrak A^{\rs}(F_{v})$ et à $G^{\rs}(F_{v})$ où $G^{\rs}$ est l'image réciproque 
de $\mathfrak A^{\rs}$.

\begin{lemme} \label{mesure nulle}
Soit $f_v$ une fonction lisse à support compact sur $G(F_{v})$, 
$dg$ une mesure invariante sur $G(F_{v})$. Notons avec les mêmes notations leurs 
restrictions à $G^{\rs}(F_{v})$. Alors, l'intégrale 
	$$\int_{G^{\rm rs}(F_{v})} f_v(g_{v}) dg_{v} $$
est absolument convergente et égale à $\int_{G(F_{v})} f_v(g_{v}) dg_{v}$.
\end{lemme}

Ceci découle du fait général que l'ensemble 
des $F_{v}$-points d'un sous-schéma fermé strict a une mesure nulle.\footnote{Nous aurions
préféré ajouter une référence
ou des références.
La démonstration de ce lemme doit être bien plus simple que la dé\-mon\-stration
générale pour un ``sous-schéma fermé strict''. Il est parfois nécessaire mais néanmoins 
dangereux d'utiliser des
résultats dont on ne comprend pas la dé\-mon\-strat\-ion et nous ne voulons pas
encourager cette habitude. Malheureusement faute de temps, nous avons accepté
de bâcler 
dans cet article plusieurs points de moindre importance.}

Soit $\gamma\in G^{\rs}(F^{\sep})$ d'image $a\in \Af^{\rs}(F^{\sep})$
et soit $T=G_\gamma$ son centralisateur. On a une suite exacte d'espaces 
tangents 
$$0 \to\rm{Tan}_{\gamma}(\mathbf c^{-1}(a)) \to \rm{Tan}_{\gamma} G \to \rm{Tan}_{b} (\Af) \to 0$$
qui induit une égalité
\begin{equation}\label{measures}
\wedge^{d} \mathfrakg=\wedge^{n} {\rm{Tan}}_{a} (\Af) \otimes \wedge^{d-n} 
\rm{Tan}_{\gamma}(\mathbf c^{-1}(a))
\end{equation}
où $d=\dim(G)$, $n=\dim(\Af)$. Ainsi le choix d'une $d$-forme invariante  $\omega_{G}$ sur $G$ 
et une $n$-forme partout non nulle $\omega_{\Af}$ sur $\Af$ induit une $(d-n)$-forme
$\omega_a=\omega_{G} \otimes \omega_{\Af}^{-1}$ sur la fibre $\mathbf c^{-1}(a)$, non nulle et 
$G$-invariante. 

Cette fibre, ou plutôt l'ensemble des points sur la fibre
à coefficients dans $F_v$, se présente de deux façons. C'est d'abord une réunion finie d'espaces
$T(F_v)\backslash G(F_v)$ donnés par $\{g^{-1}\gamma'g\}$ où $\gamma'$ parcourt un ensemble
de représentants des classes de conjugaison dans la classe de conjugaison stable
rattachée à  $\gamma$, et deuxièmement l'image inverse de $a$ par rapport à $\mathbf c$. 
La mesure traditionnelle dans sa première forme est donnée comme dans la formule \eqref{3.1}
par le quotient de mesures locales normalisées $d\bar g_v=d g_v/dt_v$. Par contre, sur l'image
inverse de $a$, il y a soit localement soit globalement la forme $\omega_a$ et
les mesures $|\omega_a|_v$
rattachées à elle.
Ces deux mesures ne sont pas les mêmes. Nous avons préféré utiliser la deuxième, qui est plus
géométrique, dans la démonstration de la proposition \ref{prodinf}. Mais c'est
la première qui est traditionnelle et aussi plus utile dans l'analyse harmonique,
donc dans la théorie des intégrales orbitales, créée par Harish-Chandra mais
avec des contributions importantes de Shalika et, dans le cadre de l'endoscopie, 
de Shelstad. Puisque nous aurons
besoin de cette théorie par la suite, il nous faudra comprendre la relation entre les
deux mesures, qui est assez simple. 

Puisque les deux mesures sont invariantes par rapport au centre $Z$, il suffit
de traiter le cas que $G=G_{\der}$, qui est selon nos hypothèses simplement connexe.
Soit $r=r_{\der}$ le rang de $G_{\der}$ et soient $\xi_1,\dots,\xi_r$ les poids
dominants des représentations $\rho_i$ de \ref{quotient adjoint}. Soit
$\omega_T=d\xi_1\wedge \dots\wedge d\xi_r$ et soit $\omega_{T\backslash G}$
une forme complémentaire invariante à droite. Elle définit alors les mesures
$d\bar g_v$. On pose
$\omega_G=\omega_T\wedge\omega_{T\backslash G}$. Nous choisissons
pour $\omega_{\mathfrak A}$, qui est sous notre hypothèse une forme sur $\mathfrak B$, à savoir
\begin{equation}\label{omega}
\omega_{\mathfrak A}=\omega_{\mathfrak B}=db_1\wedge\dots\wedge db_r,
\end{equation}
les $b_i$ étant par abus de notation des fonctions sur $T$ et sur $\mathfrak B$.

L'application $T=G_x\to\mathfrak A$ est étale au-dessus de $\mathfrak A^{\rm rs}$. L'image
inverse de $\omega_{\mathfrak A}$ est évidemment $db_1\wedge\dots\wedge db_r$.
Fixons un ordre
sur les caractères de $G_x$.
\begin{proposition}\label{HC1}
Si on choisit un ordre sur les racines et si on pose
$$
\Delta(t)=\pm t^{-\rho}\prod_{\xi>0}(\xi(t)-1),
$$
alors
$$
\omega_{\mathfrak A}=db_1\wedge\dots\wedge db_r=\pm \Delta(t)d\xi_1\wedge\dots\wedge d\xi_r=
\pm \Delta(t)\omega_T.
$$
\end{proposition}
Le signe dans cette équation n'a aucune
importance. Le carré de $\Delta(t)$ est $D(t)=\prod_\xi(\xi(t)-1)$, ce qui est une fonction
invariante et par conséquent une fonction sur $\mathfrak A$ et $\mathfrak B$. Nous écrivons
aussi $D(a)$. La fonction $\Delta(a)$ n'est définie qu'à un signe près mais nous utilisons
les normes $|\Delta(a)|_v$ qui sont bien définies.

Avant de démontrer la proposition, calculons $\omega_a$. À cette fin il faut
fixer $\omega_G$ et $\omega_{\mathfrak A}$. Cela fixe des mesures locales d'une façon
arbitraire, ce qui n'est pas grave, car il n'y a aucun effet global. Les mesures $\omega_{\mathfrak A}$
et $\omega_T$ sont données.
l'application $(\gamma,g)\to g^{-1}\gamma g$ de $T\times T\backslash G$ dans $G$
donne l'application ${\rm Tan}\times{\rm Tan}_{T\backslash G}$. Soient, avec une notation
dont l'interprétation est évidente,
$$
\omega_{T\backslash G}=\wedge_{\xi}d\xi,\qquad \omega_G=\omega_T\wedge\omega_{T\backslash G},
$$ 
où $\xi$ parcourt les racines de $G$. 
Alors l'application $t\times \bar g\mapsto g^{-1}tg=t(t^{-1}g^{-1}tg)$ de $T\times T\backslash G$
vers $G$ est étale sur $T^{\rm rs}\times T\backslash G$ et l'image de $\omega_T\wedge \omega_{T\backslash G}$
est 
$$
\{\prod_{\xi}(\xi(t)-1)\}\omega_T\wedge\omega_{T\backslash G}=
\Delta^2(t)\omega_T\wedge\omega_{T\backslash G}=\Delta^2(t)\omega_G.
$$ 
La proposition donne par conséquent 
\begin{equation}\label{conversion}
\omega_a=\pm\Delta(t)\omega_{T\backslash G}.
\end{equation}
De cette équation découle 
la relation 
\begin{equation}\label{orbit}
\int_{\mathbf c^{-1}(a)} f_v(g)|\omega_a|_v=
|\Delta(t)|_vL_v(1,\sigma_{T/G}){\rm Orb}(t_{\st},f_v), 
\end{equation}
où $\mathbf c(t)=a$.
L'intégrale à gauche a des avantages du point de vue géométrique, mais c'est l'expression à droite
qu'on emploie traditionnellement dans l'analyse harmonique locale. L'équation \eqref{orbit} est valable
si les mesures sont choisies de la façon expliquée, ce que nous supposons par la suite. 
Donc à gauche, la mesure est la mesure géométrique, tandis qu'à droite
la mesure est la mesure normalisée locale. On se perd facilement entre les deux! Les
spécialistes de l'analyse harmonique non abélienne sont habitués à l'expression de droite; les
géomètres, à celle de gauche. L'égalité
reste valable pour un groupe $G$ avec $G_{\der}$ simplement connexe. Nous posons
\begin{equation}\label{theta}
\theta_v(a;s)=L_v(s,\sigma_{T/G})|\Delta(t)|_v{\rm Orb}(t_{\st},f_v), \quad s\geq1,\quad \mathbf c(t)=a,
\end{equation} 
en observant que pour $t=\gamma\in T(F)$ régulier, l'expression \eqref{3.9} est égale à
$\theta(a;s)=\prod_v\theta_v(a;s)$, $\gamma\mapsto a$, car $\prod_v|\Delta(\gamma)|=1$. Si $\gamma$
n'est pas régulier elle est $0$.

Observons que le comportement local de $\theta_v(a;s)$, et par
conséquent le comportement global de sa transformée de Fourier,
dont l'étude sera commencée dans \cite{L5},
est fortement influencé par les facteurs $L_v(s,\sigma_{T/G})$,
car la valeur de ces fonctions dépend de la classe de
$T$, qui à son tour dépend de $a$. 

Nous nous permettons une démonstration de la proposition sur le corps de nombres complexes.
La proposition sur un corps arbitraire s'en déduit facilement. Que le signe dans la proposition \ref{HC1}
soit arbitraire est clair. Il dépend du choix de l'ordre sur les racines. 
Il est aussi
évident que $\Delta$ est une somme à coefficients entiers
de caract\`eres $t^\lambda$ de $T$ car chaque $a_i$ l'est. Ces entiers
sont indépendants du corps $F$. Il suffit donc de vérifier
la formule sur $\mathbb C$. Alors $\Delta(t)$ est \'evidemment
une fonction de $t$ alternante. Si $s$ est une r\'eflexion du
groupe de Weyl alors
$\Delta(s(t))=-\Delta(t)$, car les $a_i$ sont des fonctions
sym\'etriques et le d\'eterminant de $s$, comme application
lin\'eaire de l'espace des caract\`eres, est \'egal \`a $-1$. Il en résulte que
$\Delta(t)$ s'annule sur les vari\'eti\'es $t^\alpha=1$ dans $T$. 
On a par cons\'equent
\begin{equation}\label{alternating}
\Delta(t)=\pm Q(t)t^{-\rho}\prod_{\alpha>0}(1-t^\alpha),
\end{equation}
o\`u $Q(t)$ est une combinaison lin\'eaire des fonctions $t^\lambda$,
$\lambda$ un caract\`ere de $T$. En plus, si $s$ est la réflexion
rattachée à une racine simple $\beta$, alors 
$$
s\left(\prod_{\alpha>0,\,\alpha\neq\beta}(1-t^\alpha)\right)=
\prod_{\alpha>0,\,\alpha\neq\beta}(1-t^\alpha),
$$
et 
$$
s(t^\beta)=t^{-\beta},\qquad
s(t^{-\rho})=t^{\beta-\rho},
$$
de sorte que $Q(t)$ est une fonction invariante sous le groupe de Weyl.
Donc
\begin{equation}\label{forgotten}
Q(t)=\sum_{\lambda\geq0}a_\lambda S_\lambda(t),
\end{equation}
où $S_\lambda$ est simplement la somme sur tous les conjugués $\lambda'$
de $\lambda$ de $t^{\lambda'}$. Prenons le plus grand $\lambda$ pour lequel
$a_\lambda\neq0$. Alors le poids le plus grand à droite de
\eqref{alternating} est $\lambda+\rho$.
Puisque chaque $b_i$ est une somme semblable à \eqref{forgotten}, mais dans laquelle
le poids le plus grand pour lequel $a_\lambda\neq0$ est le poids
fondamental $\mu_i$. 

Puisque le corps de base est devenu aux fins de cette d\'emonstration
le corps des complexes, 
nous avons un mélange de notations additives et multiplicatives:
$\xi_i=e^{\mu_i(z)}$,
où $z=\sum_iz_i\hat\alpha_i$ est dans l'alg\`ebre de Lie de $T$.
Alors 
$$
d\xi_1\wedge\dots\wedge d\xi_r=
\exp(\rho(z))dz_1\wedge\dots\wedge dz_r,
$$
tandis que
$$
db_i=\exp(\mu_i(z))dz_i+\dots,
$$
où les termes qui manquent sont tous de la forme $\exp(\mu(z))$,
$\mu<\mu_i$. Il résulte alors que le terme
$\exp(\lambda(z))$ dans $\Delta(t)$ avec $\lambda$ maximal est 
$\lambda=\sum_i\mu_i=\rho$. Par cons\'equent, $Q(t)$ est une constante
et cette constante est $1$.
\vskip1pc
\noindent

Dans une exposition plus sytématique il faudrait donner une dé\-mon\-stration
valable pour tout corps de base $F$. Nous ne sommes pas aussi ambitieux à 
ce stade-ci.

\section{Adélisation de la formule des traces}
Les fonctions $\theta_v(a;s)$ de la formule \eqref{theta}
sont des fonctions sur $\mathfrak A_v$ à support compact.
Il est bien connu qu'elles sont bornées mais pas nécessairement lisses.
Nous offrons comme références, celles que nous utiliserons dans des articles
à suivre: pour les groupes 
sur $\mathbb R$ et $\mathbb C$ le livre \cite{V} de Varadarajan; pour les corps
$p$-adiques l'article \cite{Sh} de Shalika;
pour les groupes sur les corps locaux de caract\'eristique positive, nous n'avons
pas de référence. En plus des questions locales,
le comportement du produit 
$$
\theta(a;s)=\prod_v\theta_v(a_v;s), \qquad a=\prod_va_v,
$$ 
n'est pas
\'evident m\^eme pour $s=1$ et $f_v$ \'egal \`a la fonction caract\'er\-istique
de $G(\calO_v)$. En particulier, la transform\'ee de Fourier locale,
$\hat\theta_v(a;s)$, n'est pas \`a support compact et la convergence de l'intégrale
qui définit la transform\'ee
de Fourier globale $\hat\theta(a;s)=\prod_v\hat\theta_v(a;s)$ n'est pas acquise. Puisque 
nous proposons d'utiliser
la formule de Poisson, le comportement de $\theta_v$ et 
de $\theta$ exige une \'etude plus pouss\'ee. Nous la commen\c cons
dans cet article et pour
le groupe $SL(2)$ on la poursuivra dans un prochain article. On espère
revenir \`a la question g\'en\'erale assez tôt.  

La formule de Poisson dans sa forme habituelle ad\'elique est une application de l'analyse
harmonique \`a la paire autoduale $F\subset\bbA$. Le mauvais comportement de la fonction
sugg\`ere l'utilisation d'une forme tronqu\'ee, et plus pr\'ecisement dans une forme sugg\'er\'ee
par une observation de Jayce Getz, qui avait dans \cite{G} d\'ej\`a
rencontr\'e une difficult\'e semblable. 
\begin{lemme}\label{getz}
Soit $T$ un tore sur le corps global $F$. Soit $\mathfrak X$ le groupe des caract\`eres
rationnels de $T$, un groupe pourvu d'une action du groupe de Galois $\Gal(F^{\sep}/F)$. Supposons que
$\Lambda=\{\lambda_1,\dots,\lambda_n\}\subset\mathfrak X$ et que $\Lambda$ soit
invariant par rapport \`a $\Gal(F^{\sep}/F)$. Supposons enfin que $S$ soit un ensemble fini
de places de $F$ qui contient toutes les places infinies et que, pour chaque $v\in S$, 
un
sous-ensemble compact $U_v\subset T(F_v)$ soit donné. Alors il existe un sous-ensemble fini $S'\supset S$
tel que si $t\in T(F^{\sep})$, $t\in U_v$ pour tout $v\in S$, et $|\lambda(t)|_v=1$
pour tout $v\notin S$ et tout $\lambda\in\Lambda$,
alors
$$
\prod_{j=1}^n\left(1-\lambda_j(t)\right)=0
$$
ou
$$
|1-\lambda_j(t)|_v=1
$$
pour tout $j$ et tout $v\notin S'$.
\end{lemme}

Quoique le lemme est important, sa d\'emonstration est facile. Il y a certainement
un nombre positif $A$ tel que $1/A\leq |\lambda(t)|_v\leq A$
pour tout $\lambda\in\Lambda$ et  tout $v\in S$. En plus $|1-\lambda_j|_v\leq 1$ pour $v\notin S$.
Consid\'erons $\alpha=\prod_j(1-\lambda_j(t))$. Il appartient \`a $F$. Pour $v\notin S$,
$|\alpha|_v\leq1$, et pour $v\in S$, $|\alpha|_v\leq A^n$. Il n'y a qu'un nombre fini
de places en dehors de $S$ telles que $F_v$ contient un \'el\'ement de valeur absolue
positive mais plus grande ou \'egale \`a $A^{-n}$, donc telles que $q_v\leq A^n$. Soit $S'$
leur r\'eunion avec $S$. L'ensemble $S'$ satisfait aux conditions du lemme
et ne d\'epend que de $S$, $A$ et $n$. Il est
évident que l'ensemble $S'$ devient de plus en plus grand, et à une allure inquiétante,
lorsque $S$ croît ou
les ensembles $U_v$, $v\in S$,
grandissent.
    
Le lemme 4.1 permet une nouvelle formulation de la formule des traces qui nous sera importante.
Nous commen\c cons avec une observation \'el\'ementaire
sur la formule de Poisson ad\'elique. Soit $S'$ un ensemble
fini de places de $F$ qui contient toutes les places infinies. Fixons un caract\`ere global $\chi$
de $\bbA_F$ avec la propri\'et\'e habituelle:  l'ensemble des $b\in\bbA_F$
tels que $\chi(ba)=1$ pour tout $a\in F$ est $F$. Nous supposons que $S'$ soit suffisamment grand
pour que la mesure auto-duale sur $F_v$, $v\notin S'$, donne \`a $\calO_v$ la mesure $1$.
Donc l'ensemble des $x\in F_v$ tels que $\chi_v(x)=1$ est $\calO_v$. Soit $\bbA_F^{S'}$ le sous-ensemble
des $a\in\bbA_F$ tels que $a_v\in\calO_v$ pour $v\notin S'$. Soient 
$$
\calO^{S'}=\prod_{v\notin S'}\calO_v,
\qquad\mathbb A_{S'}=\prod_{v\in S'}F_v 
\qquad F_{S'}=F\cap \bbA_F^{S'}. 
$$
Nous identifions $F_{S'}$
avec son image dans $\mathbb A_{S'}$. 
\begin{lemme}\label{getz1}
Si $S'$ est suffisamment grand alors $\calO_v$ est auto-dual par rapport à $\chi_v$ pour $v\notin S'$
et
\begin{equation}\label{getz2}
\bbA_F=\bbA_F^{S'}+F.
\end{equation}
En plus $F_{S'}\subset\bbA_{S'}$ est auto-dual par rapport \`a la restriction de $\chi^{S'}$
\`a $\bbA_{S'}$ et la mesure produit sur $\bbA_{S'}$ donne 
${\rm meas}(F_{S'}\backslash\bbA_{S'})=1$.
\end{lemme}
Il est évident que ce lemme implique un lemme semblable pour n'importe quel espace
vectoriel de dimension finie sur $F$. Nous l'utiliserons dans cette forme plus générale.

Nous disons que $S'$ est suffisamment grand s'il contient un sous-ensemble fini convenable
que nous
allons construire. Pour simplifier la notation nous dénotons ce
sous-ensemble $S'$. Il est évident que l'équat\-ion \eqref{getz2} est la clé
car elle implique que
$$
\bbA_F=F+\calO^{S'}+\bbA_{S'}. 
$$
Un caractère de $\bbA_{S'}$ trivial sur $F_{S'}$ s'étend d'une façon unique à un caractère
de $\calO^{S'}+\bbA_{S'}$ trivial sur $\calO^{S'}+F_{S'}$ et ensuite à un caractère
de $\bbA_F$ trivial sur $F$.
Il suffit donc
de d\'emontrer le lemme pour un seul $S'$. Il sera alors valable pour tout ensemble plus
grand. 

Le corps $F$ est une extension finie et s\'eparable d'un corps $F^0$, soit le corps des nombres
rationnels, soit le corps de fonctions sur la droite projective $\bbP^1$ sur un corps fini.
Nous pouvons choisir pour $S'$ l'ensemble des places au-dessus d'un ensemble $S_0'$
de places de $F^0$. Nous choisissons cet ensemble tel qu'il existe une base $\{\alpha_1,\dots,\alpha_d\}$
de $F/F^0$ pour laquelle le d\'eterminant ${\rm det}(\tr\alpha_i\alpha_j)$ soit de valeur absolue $1$ pour tout
$v$ en-dehors de $S_0$ et pour laquelle chaque $\alpha_i$
appartienne à $F_{S'}$. Il est alors \'evident que
$$
F_{S'}=\{\sum_ia_i\alpha_i\,\vert\,a_i\in F_{S'_0}^0\},\qquad
\bbA_F^{S'}=\{\sum_ia_i\alpha_i\,\vert\,a_i\in \bbA_{F^0}^{S'_0}\}.
$$
Il suffit par cons\'equent de d\'emontrer le lemme pour $F=F^0$. Pour $F=\QQ$, il est
\'evident que le lemme est valable si $S'=\{\infty\}$. Pour le corps de fonctions de $\bbP^1$,
on peut prendre pour $S'$ une seule place, n'importe laquelle.

La partie principale de la formule des traces est une somme sur les classes de conjugaison
elliptiques et régulières. Nous avons observé qu'il y a un facteur commun à tout terme de la somme,
le facteur $m_G\rho_G\tau(G)$ qui est égal à $1$ si $G$ est semi-simple et simplement connexe.
Pour simplifier les formules
nous n'incluons pas ce facteur par la suite. En effet, puisque le but principal de cet article
est d'introduire la somme donnée par la formule des traces comme une somme de Poisson, il
est convenable de ne considérer désormais que les groupes semi-simples et simplement connexes. En général
on a une somme sur $\eta\in\mathfrak h$, sur $\mathfrak B_\eta(F)$ et sur $Z_\eta(F)$. C'est la somme
sur $\mathfrak B_\eta(F)$ qui pose des difficultés et pour laquelle on utilise la dualité de Poisson.
Les paramètres $z\in Z_\eta(F)$ et $\eta$ ne sont à toutes fins utiles que des indices. Il est
donc préférable de les écarter en supposant que $G=G_{\der}$ est semi-simple et simplement
connexe.

Grâce au lemme \ref{getz}, pour une fonction $f=\prod f_v$, où les $f_v$ sont lisses et à support compact,
le produit 
$$
\prod_v{\rm Orb}(\gamma_{\st},f_v)
$$
est $0$ sauf pour un nombre fini de classes de conjugaison stables et régulières. Dans cet article
comme il arrive souvent avec la formule des traces, les classes singulières sont mises à part pour 
être traitées plus tard lorsque les grandes lignes de l'argument sont plus claires. Elles ne le sont
pas encore. 
Donc
$$
\sum_\gamma L(1,\sigma_{T/G})\prod_v{\rm Orb}(\gamma_{\st},f_v)=\lim\limits_{s\searrow 1}\sum
L(s,\sigma_{T/G})\prod_v{\rm Orb}(\gamma_{\st},f_v).
$$
Puisque $|\Delta(\gamma)|=\prod_v|\Delta(\gamma)|_v=1$ pour les classes semi-simples
et régulières, cette expression est égale à
\begin{equation}\label{prod1}
\lim\limits_{s\searrow 1}\sum\theta(\mathbf c(\gamma_{\st});s).
\end{equation}
Dans la somme, $a=\mathbf c(\gamma_{\st})$ ne parcourt encore 
qu'un sous-ensemble de l'ensemble $\mathfrak B(F)$, l'indice $\eta$ étant omis
maintenant qu'il n'y en a qu'un seul. Ce qui manque, ce sont les $\gamma_{\st}$
qui ne sont pas elliptiques.

Choisissons un ensemble $S'$ qui satisfait aux conditions des lemmes \ref{getz} et \ref{getz1}.
Nous supposons aussi que, en dehors de $S'$, le groupe $G$ est quasi-déployé et déployé sur une extension
non ramifiée et que $f_v$ est la fonction caractéristique d'un sous-groupe hyperspécial.
Les conditions du deuxième de ces lemmes sont indépendantes des fonctions $f_v$
utilisées dans la formule; celles du premier dépendent du support de ces fonctions.
Lorsque nous commençons l'analyse des sommes de Poisson il faudra en tenir compte.
Dans cet article ces fonctions seront fixées aussi bien que $S'$, mais il est certain
que le fait que $S'$ dépende des fonctions $f_v$ rendra ultérieurement l'analyse plus difficile.

Nous avons choisi $S'$ tel que si $v\notin S'$ et si l'intégrale orbitale ${\rm Orb}(\gamma_{\st},f_v)\neq0$
alors $|\Delta(\gamma)|_v=1$. Nous rappelons quelques formules des articles \cite{O1,O2,O3}. Soit
$\kappa_v$ le corps résiduel par rapport à $v$.
D'abord pour un tore $T$ avec bonne réduction en une place $v$ où $\mathfrak d_v=\calO_v$,
le côté droit de la formule \eqref{Weil} est égal
à
\begin{equation}\label{tore}
q_v^{-\dim T}|T(\kappa_v)|=L(1,\sigma_T)^{-1},
\end{equation}
cette égalité étant la formule (1.2.6) de \cite{O1}. 
Nous vérifions en plus que sous nos hypothèses sur $\gamma$, c'est-\`a-dire,
\vskip.1pc
\noindent
(i) $|\lambda(\gamma)|_v=1$ pour tout caractère $\lambda$ du sous-groupe
de Cartan $T$ et
\vskip .1pc
\noindent
(ii) $|1-\alpha(\gamma)|_v=1$
pour toute racine $\alpha$, 
\vskip.1pc
\noindent
on a
\begin{equation}\label{genorb}
{\rm Orb}(\gamma_{\st},f_v)=q_v^{-\dim G_{\der}}|G_{\der}(\kappa_v)|.
\end{equation}

Cette équation résulte de deux conséquences des hypothèses (i) et (ii) sur $\gamma$: 
\vskip.1pc
\noindent
(a) l'élément $\gamma$ lui-même
n'est contenu que dans un seul sous-groupe hyperspécial;
\vskip.1pc
\noindent 
(b) 
la classe de conjugaison stable de $\gamma$ a une
intersection non vide avec
tout sous-groupe hyperspécial $K$ de $G_{\der}$. 
\vskip.1pc
\noindent
Ceci est une conséquence des propriétés des
immeubles rattachés à $G(F_v)$. Malheureusement, mais certainement
faute d'efforts, en dépit de l'usage presqu'universel
des immeubles, nous n'avons pas trouvé de références qui donnent exactement ce que nous
cherchons. Nous continuons néanmoins. 
Supposons que la fonction $f_v$ soit la fonction caractéristique du sous-groupe $K$ et que $K$
soit hyperspécial. Soit $x$ le point de l'immeuble de Bruhat-Tits qui définit $K$.
En remplaçant $\gamma$ par un élément auquel il est stablement conjugué,
nous pouvons supposer qu'il est contenu dans $K$. Si $\gamma'$ est stablement conjugué
à $\gamma$ alors $\gamma'=g^{-1}\gamma g$, $G\in G(F^{\rm sep})$. Si $\gamma'$ fixe $K$ ou $x$, alors
il résulte de la première des trois conséquences, mais pour une extension finie de $F$, que $g$ fixe
$x$, mais le cocycle $\sigma\rightarrow g^{\sigma-1}$ prend alors ses valeurs dans $K=G(\calO_v)$.
Il résulte du théorème de Lang qu'il est trivial et que $\gamma'$ et $\gamma$ sont conjugués dans $K$.
Par conséquent,   
\begin{align*}
{\rm Orb}(\gamma_{\st},f_v)&=\frac{1}{L_v(1,\sigma_{T/G})}\frac{{\rm meas}_{\geom}K_v}{\meas_{\geom}(K_v\cap T)}\\
&=q_v^{-\dim G_{\der}}|G_{\der}(\kappa_v)|.
\end{align*}
Rappelons (\cite{O3}) que
le produit
$$
\Pi_{S'}=\prod_{v\notin S'}q_v^{-\dim G_{\der}}|G_{\der}(\kappa_v)|
$$
converge. Il est évidemment indépendant de $T$. Nous posons
$$
\Pi_{S'}(\gamma_{\st},sd)=L_{S'}(s,\sigma_{T/G})\Pi_{S'}.
$$

Considérons la formule des traces, ou plutôt sa partie
elliptique, qui est devenue après nos transformations
la formule \eqref{prod1}, pour 
une fonction $f=\prod f_v$ donnée et choisissons
l'ensemble $S'$ suffisamment grand. Il résulte du choix de $S'$ que $|\Delta_v(\gamma)=1|$
pour $v\notin S'$ de sorte que
\begin{equation}\label{theta1}
\theta(\gamma_{\st};s)=L_{S'}(s,\sigma_{T/G})\Pi_{S'}\prod_{v\in S'}\theta_v(\gamma_{\st};s).
\end{equation}
Nous avons introduit la troncation de Getz pour éviter une seule difficulté: nous ne savons
pas vérifier que le comportement de la transformée de Fourier de $\theta(a;s)$, $s>1$,
permet l'application de la dualité de Poisson. Une telle possibilité est peu
probable. L'équation \eqref{theta1} soulève une autre 
difficulté car le facteur $L(s,\sigma_{T/G})$ fait apparaître dans \eqref{theta1} le facteur
$L_{S'}(s,\sigma_{T/G})$ qui dépend de $\gamma$ par l'intermédiaire de $T$. Puisque pour
un $f$ donné il n'y a qu'un nombre fini de $\gamma$ nous pouvons, en choisissant $S'$
suffisamment grand, le choisir tel que le  
facteur $L_{S'}(s,\sigma_{T/G})$ est sur un intervalle ouvert $(1,1+\epsilon)$ aussi proche de $1$
que désiré et cela sans modifier le dernier facteur de \eqref{theta1}. Le deuxième facteur change
d'une façon uniforme, donc indépendamment de $\gamma$, et il s'approche de $1$ lorsque
$S$ grandit. Puisque l'ensemble des $\gamma$ qui apparaissent dans \eqref{3.9} ne change pas
lorsqu'on fait grandir $S'$, nous pouvons remplacer la somme des termes \eqref{3.1} 
par la somme sur le même ensemble des $\gamma$ de
\begin{equation}\label{theta2}
\Pi_{S'}\prod_{v\in S'}\theta_v(\gamma_{\st};1)=
\Pi_{S'}\lim_{s\searrow 1}\prod_{v\in S'}\theta_v(\gamma_{\st};s).
\end{equation}
Si on tient compte du fait que $\Pi_{S'}$ s'approche de $1$ nous pouvons aussi et pour les
mêmes raisons le supprimer.\footnote{Nous offrons à ce stade-ci une explication brève de notre
usage du mot {\it adélisation}. Les espaces vectoriels $V$ sur le corps $F$ et leur
produit tensoriel $\bbA_F\times V$ interviennent dans la théorie analytique des
nombres pour étudier, en utilisant la formule
de Poisson, les sommes $\sum_{v\in V}\varphi(v)$, $\varphi=\prod_v \varphi_v$ et leur comportement
asymptotique. Cette formule est un des seuls outils disponibles. Nous voulons 
l'utiliser pour $\theta(a;s)$.
Cela nous est interdit car les fonctions $\theta_v$ ne sont pas lisses. Il apparaît
par contre que leur comportement est assez bon pour que leur utilisation soit permise
non pas dans le cas des produits infinis restreints des corps locaux mais 
dans le cas des produits tronqués.
Malheureusement nous ne sommes pas en état de tout expliquer en quelques pages. Pour être
francs, nous n'avons guère nous-mêmes entamé l'étude des problèmes qui ressortent de ces questions
(\cite{L5}).}

Nous avons escamoté toutefois une petite difficulté. C'est que le produit infini qui définit
$L_{S'}(s,\sigma_{T/G})$ ne converge pas uniformement sur l'intervalle. Ce qui converge
uniformément est le produit extérieur dans
$$
\prod_{n=1}^\infty\left\{\prod_{n<q_v\leq n+1}L_v(s,\sigma_{T/G})\right\}.
$$
Il s'agit d'une forme du théorème dit de Chebotarev dont une forme se
trouve dans \cite{IK}. Nous n'avons pas cherché d'autre référence. Il y en a
beaucoup. 
Nous supposons par conséquent que $S'$ est toujours un ensemble défini comme
$\{v\,|\,q_v\leq n\}$, ce qui ne l'empêche pas d'être suffisamment grand. 

Pour étudier la somme des termes de \eqref{theta2}
qui apparaît dans la formule des traces, mais avec $\Pi_{S'}$
supprimé, à savoir,
\begin{equation}\label{poisson}
\sum \prod_{v\in S'}\theta_v(\gamma_{\st};s),
\end{equation}
nous pouvons utiliser la dualité de Poisson pour la
paire $F_{S'}\subset \bbA_{S'}$. En effet, selon le lemme \ref{getz} et le choix de $S'$,
$\mathbf c(\gamma_{st})$ est entier en dehors de $S'$ si $\theta_v(\gamma_{\st},1)\neq0$
pour tout $v\in S'$. Par conséquent, la somme \eqref{poisson} s'écrit
\begin{equation}\label{missing}
\sum_{b\in \mathfrak B_{S'}}\prod_{v\in S'}\theta_v(b;s),
\end{equation}
où $\mathfrak B_{S'}$ est défini comme $\Spec F_{S'}[b_1,\dots,b_r]$.
 
Notre affirmation que nous pouvons utiliser la dualité de Poisson est optimiste
et, en plus, pas encore tout à fait correcte.
D'abord les fonctions $\theta_v(b;s)$ ne sont pas lisses. Elles sont singulières
sur la sous-variété $D(b)=0$. Il est toutefois
probable qu'elles sont suffisament lisses, donc que leurs singularités sont assez
modérées, pour que la somme de Poisson duale, 
\begin{equation}\label{hatpoisson}
\sum_{b\in\mathfrak B_{S'}}\prod_{v\in S'}\hat\theta_v(b;s),
\end{equation}
converge absolument et soit égale à la somme de Poisson pour la fonction
$\prod_{S'}\theta_v(b;s)$.
L'un de nous espère vérifier ceci pour $SL(2)$ dans un prochain article (\cite{L5}). Il n'y a pas de
raison de croire que le résultat sera tellement difficile même pour un groupe arbitraire.
Les fondements de la théorie nécessaire se trouvent en toute probabilité dans \cite{{Sh}, V}.

Il sera utile par la
suite d'avoir des notations plus précises,
\begin{equation}\label{precise}
\theta_{S'}(b;s)=\prod_{v\in S'}\theta_v(b;s),\qquad\hat\theta_{S'}(b;s)=\prod_{v\in S'}\hat\theta_v(b;s).
\end{equation}

Une deuxième difficulté est que nous avons écarté dès le début les $b$
qui ne définissent pas un tore elliptique. Nous pouvons les ajouter à \eqref{missing}.
Mais alors il faut ensuite les soustraire,
car ceux-ci posent des problèmes de convergence. En particulier pour eux la fonction $L_{S'}(s,\sigma_{T/G})$
converge vers l'infini lorsque $S'$ devient de plus en plus grand et
$s\searrow1$. On peut toutefois supposer que ces divers infinis
se compenseront. Par exemple pour le groupe $SL(2)$, il est naturel en examinant le comportement des
sommes telles que \eqref{RL.1.1} de soustraire la contribution des représentations qui sont
les transferts de représentations des groupes diédraux et d'y ajouter même le tore dé\-ployé, ce qui
rend le processus d'addition et de soustraction plus naturel. Il reste néanmoins un casse-tête.
Comment dans la formule ainsi obtenue s'annulent mutuellement les deux infinis. On revient
à cette question dans \cite{L5}. 

Pour le faire il faudra introduire 
et étudier la notion
de transfert stable $f=f^G\mapsto f^H$ et des facteurs de transfert. À la fin on étudiera
non pas les sommes de Poisson pour la seule fonction $\theta$ ou $\hat\theta$ mais
pour des différences $\theta-\varphi$ et $\hat\theta-\hat\varphi$, où $\varphi$
est défini à partir des transferts $f^H$. Même sans $\varphi$ nous pouvons examiner
le terme $\hat\theta(0;s)$ qui reste fini lorsque $s\rightarrow1$ même si on ajoute
les termes rattachés aux tores non elliptiques. Nous verrons, mais non pas dans cet
article, que la contribution de $\hat\varphi(0;s)$ sera $0$.

Avant de commencer rappelons la suite de modifications. Nous introduisons le paramètre $s>1$,
qu'il faut faire converger vers $1$. Ensuite nous introduisons l'ensemble $S'$ suffisamment grand
et nous passons à la somme tronquée. Pour cette somme tronquée nous utilisons la dualité de Poisson.
Mais alors le terme 
$$
\hat\theta_{S'}(0)=\lim_{s\rightarrow1}\hat\theta_{S'}(0;s)
$$ 
dans la somme duale dépend de $S'$. Pour obtenir
une valeur in\-dé\-pend\-ante de $S'$, il faut passer à la limite double\footnote{Un rapporteur
a été troublé par cette limite double. Avec raison. En effet l'ordre des limites est indifférent.
Ces limites, en particulier, la limite $s\searrow$ sont discutées avec plus de soin dans \cite{L5}} 
\begin{equation}\label{double}
\hat\theta(0)=\lim\limits_{\substack{
S'\rightarrow\infty\\s\searrow 1}}\hat\theta_{S'}(0;s),
\end{equation}
dans laquelle $S'$ devient de plus en plus grand de la façon prescrite.
C'est cette valeur limite dont il est question dans la dernière section de cet article.
Il faudra passer à cette limite double non pas seulement pour le terme principal,
mais aussi pour la somme $\sum_{b\neq0}\hat\theta_{S'}(b;s)$ qui reste. Pour
elle il faudra traiter non pas ce terme lui-même mais
la différence
$$
\sum_{b\neq0}\hat\theta_{S'}(b;s)-\sum_{b\neq0}\hat\phi_{S'}(b;s).
$$
On n'entreprendra l'étude des termes $\hat\phi_{S'}(0;s)$ que dans un prochain article.
 
\section{Le terme dominant}

À présent l'évidence la plus persuasive de la promesse de notre stratégie
est qu'elle permet d'isoler la contribution dominante pour $G$  semi-simple
et simplement connexe;
en effet, pour n'importe quel groupe qui satisfait à nos conditions,
il n'y a que
des questions formelles qui distinguent le cas $G=G_{\der}$ du cas
général. Il s'avère
que ce terme dominant est donné par $\hat\theta(0)$.

Nous supposons que $G=G_{\der}$. En plus, $G$ est quasi-déployé
presque partout, même partout. Par conséquent le sous-groupe dérivé
$G_{\der}(F_v)$ est $G_{\der}(F_v)$
lui-même. \big(Nous n'avons pas cherché la meilleure référence, mais
\cite{St2} est une possibilité.\big) Il en résulte que
la seule représentation automorphe de dimension $1$ de $G_{\der}(\bbA_F)$
est la représentation triviale $\pi_0$. Selon la formule \eqref{1.2}, le terme
dominant pour une représentation $\rho$ de ${}^LG$ et pour la fonction
$f$ de \eqref{f} est alors
\begin{equation}\label{zeta}
\prod_i\zeta_S(s+i)\prod_{v\in S}\tr\pi_v(f_v)
\end{equation}
où la variable $s$ est comme dans la section 1. Les entiers $i$ sont
ceux donnés par le plongement $\phi$
de $SL(2)$ dans ${}^LG$ rattaché à l'élément unipotent
principal. Les $i$ sont la moitié des poids de
$$
z\rightarrow \rho\cdot\phi(\left(\begin{matrix}
z&0\\0&z^{-1}\end{matrix}\right)).
$$
L'expression \eqref{zeta} est égale à $\tr(\pi_0(f))$,
donc à
$$
\int_{G(\bbA_F)}f(g)dg=\prod_v\int_{G(F_v)}f_v(g_v)dg_v.
$$
La prochaine formule découle de la formule \eqref{orbit},
\begin{equation}
\int_{G(F_v)}f_v(g_v)dg_v=\int_{\mathfrak A(F_v)}\theta_v(b_v;1)db_v,
\end{equation}
\label{formula1}
ou plus généralement de
\begin{equation}\label{formula2}
\int_{G(F_v)}\frac{L_v(s,\sigma_{T/G})}{L_v(1,\sigma_{T/G})}f_v(g_v)dg_v=\int\theta_v(b_v;s)db_v=\hat\theta_v(0;s).
\end{equation}
Nous soulignons que la fonction dans cette intégrale dépend de $b_v$ ou de $g_v$
par l'intérmédiaire de $T=T_{g_v}$, le centralisateur de l'élément régulier $g_v$. Il
est convenable de l'écrire comme
$$
f_{v}(g_{v};s)= \frac{L_v(s,\sigma_{T/G})} 
{L_v(1,\sigma_{T/G})} f_v(g_v). 
$$
Observons que le quotient $L(s,\sigma_{T/G})/L(1,\sigma_{T/G})$
est borné pour $1\leq s\leq 1+\epsilon$ car il est le produit
de facteurs $(1-\alpha/q_v)/(1-\alpha/q_v^s)$ où $|\alpha|=1$ où les $\alpha$ appartiennent
à un ensemble fini de racines de l'unité. Le seul cas même un peu délicat
est $\alpha=1$ mais alors $1-1/q_v\leq 1-1/q_v^s$ pour $s\geq1$.

Le lemme suivant est équivalent à la formule \eqref{genorb} mais dans
une forme qui
plaît plus aux géomètres.
\begin{lemme}
Soit $v\in |F|$ une place où $G$ a une réduction réductive. Le
morphisme ${\bf c}:G \to \Af$
s'étend alors à $\mathcal O_{v}$. Supposons que les formes volumes
$\omega_{G}$ et $\omega_{\Af}$
s'étendent en des
formes partout non nulles sur leur $\mathcal O_{v}$-modèle. Soit
$f_{v}$ la fonction caractéristique
du compact $G(\Oc_{v})$. Pour tout $b_{v}\in \Af^{\rs}(\Oc_{v})$, on a alors
\begin{equation}
       \theta_{v}(b_{v};1)=q_v^{-\dim(G)+\dim(T)} \frac{|G(k_{v})|}
       {|T(k_{v})|}=L_v(1,\sigma_{T/G})\frac{G_{\der}(\kappa_v)}{q_v^{\dim(G_{\der})}}.
\end{equation}
Ici $T$ est le centralisateur d'une
section $g_{v} \in G^{\rs}(\mathcal O_{v})$ au-dessus de $b_{v}$ et 
un tore défini sur $\mathcal O_{v}$. 
\end{lemme}

Le résultat qu'il s'agit de démontrer est le suivant.
\begin{proposition}\label{prodinf} On a l'égalité,
$$
\lim\limits_{\substack{
S'\rightarrow\infty\\s\searrow
1}}\hat\theta_{S'}(0,s)=\int_{G(\mathbb A)}f(g)dg.
$$
\end{proposition}

Le produit infini
$$
\prod_v\int_{G(F_v)}f(g_v)dg_v
$$
converge. En plus, pour un $v$ donné,
$$
\lim_{s\rightarrow1}\hat\theta_v(0,s)=\int_{G(F_v)}f_v(g_v)dg_v.
$$
Il s'agit donc de montrer que pour un $S''$ donné fini,
$$
\lim\limits_{\substack{
S'\rightarrow\infty\\s\searrow 1}}\{\prod_{
S'-S''}\hat\theta_v(0,s)-\prod_{S'-S''}\int_{G(F_v)}
f_v(g_v)dg_v\}=0.
$$
On choisit naturellement $S''$ de sorte que tous les ennuis possibles
sont écartés.

Il s'agit de deux produits, dont le deuxième vaut
$$
\prod_{S'-S''}\int_{G(F_v)} f_v(g_v)dg_v=
\prod_{S'-S''}\frac{|G_{\der}(\kappa_v)|}{q_v^{\dim
G_{\der}}}=\prod_{S'-S''}(1+O(q_v^{-2})),
$$
de sorte qu'il converge absolument.
Rappelons que nous avons supposé que $G=G_{\der}$. Pour terminer la
démonstration il
ne faut qu'une majoration semblable
\begin{equation}\label{estimate}
\hat\theta_v(0,s)=1+O(q_v^{-3/2})
\end{equation}
des facteurs dans le premier produit.

Nous avons constaté que les fonctions $f_v(g_v;s)$ sont uniformément majorées
et portées par les compacts $G(\calO_v)$. 
Toutefois,
elles sont assez compliquées et deviennent très irrégulières lorsqu'on
s'approche du diviseur donné par le carré $D$ de la fonction $\Delta$. Puisque 
cette fonction est invariante, elle
est aussi une fonction, que nous appelons le discriminant, sur la base de Steinberg-Hitchin.
Pour le groupe $SL(2)$, si le polynôme caractéristique est $X^2-bX+1$ et les valeurs
propres, $\alpha$ et $\alpha^{-1}$, alors $D=(\alpha^2-1)(\alpha^{-2}-1)=4-b^2$
est un diviseur sur la ligne droite à deux points simples, au moins si la caractéristique
n'est pas égale à $2$! Nous escamotons ce cas dans cet article.

La dé\-mon\-stration de \eqref{estimate} est fondée sur les observations suivantes. D'abord,
sur l'ouvert compact $G^{\rm tvl}(\Oc_{v})$ (l'indice ``tvl''
veut dire transversal) de
$G(\Oc_{v})$ des éléments $g_{v}$ ayant une réduction $\bar g_{v}$
régulière, mais pas
nécessairement semi-simple et dont le discriminant $D(g_v)$
a une valuation plus petite ou égale à $1$, la valeur $f_{v}(g_{v};s)$
ne dépend que de $\bar g_{v}$.
Deuxièmement, cette valeur s'exprime facilement à l'aide d'un faisceau
$\ell$-adique sur le groupe $G$ sur le
corps $\kappa_{v}$. Enfin, le complément de cet ouvert dans
$G(\Oc_{v})$ est de mesure trop petite
pour nuire à la majoration \eqref{estimate}.

Quoique les géomètres n'en ressentiront pas le besoin,
pour faciliter la com\-pré\-hension de l'argument par les spécialistes
de la formule
des traces ou des formes automorphes, nous rappelons au fur et à mesure 
les propriétés que nous
utilisons, en nous plaçant dans le cadre du groupe $SL(2)$ et en les exprimant d'une façon
concrète . Par exemple, pour $SL(2)$ l'ouvert $G^{\rm tvl}(\Oc_{v})$ 
est défini par la condition qu'une valeur propre $\alpha$ satisfait à $|1-\alpha^2|_v\leq q_v^{-1/2}$.
Donc, si la caractéristique résiduelle n'est pas $2$, les $\alpha$ qui sont exclus, ou plutôt pour lesquels
un argument plus délicat est exigé, sont ceux pour lesquels $|1\mp\alpha|_v=1,q_v^{-1/2}$.
Donc, comme en général, pour $SL(2)$ {\it transversal} veut dire que les valeurs propres
de $g$ dans $\mathfrak g_g\backslash\mathfrak g$ assument des
valeurs suffisamment distantes de $1$. Cependant pour le cas général, la définition précise
est plus compliquée.

\begin{lemme}
Le complémentaire de $G^{\rm tvl}(\Oc_{v})$ dans $G(\Oc_{v})$ est de
mesure $O(q_{v}^{-2})$.
\end{lemme}

\begin{proof}
Soit $G^{\rm sing}$ le fermé de Zariski de $G$ complémentaire
de l'ouvert $G^{\rm reg}$ des éléments réguliers, semi-simple ou pas.
D'après Steinberg (\cite{St1})
on sait que $G^{\rm sing}$ est un fermé de codimension trois. Par
exemple, pour $G=SL(2)$
il est $\{\pm I\}$.

La fonction discriminant $D$ définit un diviseur réduit
${\rm div} D$ sur $T/W$, donc sur la base de Steinberg-Hitchin $\mathfrak A$. Notons ${\rm div} D^{\rm sing}$
le sous-schéma fermé de ${\rm div}D$ où le discriminant s'annule avec un ordre au moins égal à 
deux, donc où sa différentielle est nulle. Ce fermé 
est généralement de codimension deux, 
donc vide si $G=SL(2)$, mais il y a des exceptions. Sauf pour 
ces exceptions, son image réciproque dans
$G^{\rm reg}$ est de codimension au moins deux puisque le morphisme
$\mathbf c:G^{\rm reg}\to \mathfrak A$ est lisse. Puisque $G^{\rm sing}$ est
de codimension trois, la réunion de $\mathbf c^{-1}({\rm div}D^{\rm sing}) \cup
G^{\rm sing}$ est un sous-schéma fermé de codimension au moins deux de $G$.

Le cas exceptionel le plus simple est le groupe $SL(2)$ lorsque la caractéristique
est $2$ et $D=b^2$ car alors ${\rm div} D^{\rm sing}$ contient le point $b=0$. Pour
éviter ce genre de problème il suffit d'exiger que la caractéristique
soit plus grande que $2$. L'argument qui suit est alors
valable. Nous soulignons cependant que la théorie recherchée est censée être
valable sans exception. Malheureusement pour le moment nous n'avons pas de démonstration
complète, mais nous avons déjà mis de côté pour une autre occasion
quelques cas de petite caractéristique. 

Pour tout $\bar g_{v}\in G(\kappa_{v})$, l'ensemble des
éléments $g_v\in G(\Oc_{v})$ de ré\-duct\-ion $\bar g_{v}$ est de
mesure $q^{-\dim(G)}$. Par conséquent, sous nos hypothèses, l'ensemble des $g_v\in
G(\Oc_{v})$ ayant une ré\-duct\-ion dans
$\mathbf c^{-1}({\rm div}D^{\rm sing}) \cup G^{\rm sing}$ est de
mesure $O(q_{v}^{-2})$
lorsque $q_{v}$ tend vers l'infini.

Notons $G'$ l'ouvert complémentaire de $\mathbf c^{-1}({\rm div}D^{\rm sing})
\cup G^{\rm sing}$ dans $G$.
Soit $g_{v}\in G(\Oc_{v})$ de réduction $\bar g_{v}\in
G'(\kappa_{v})$. Si $\bar g_{v}$ n'appartient pas à
$\mathbf c^{-1}({\rm div}D-{\rm div}D^{\rm sing})$, le discriminant de $\bar g_{v}$ a
la valuation nulle et $g_v$ appartient à $G^{\rm tvl}$.
Si $\bar g_{v}\in \mathbf c^{-1}({\rm div}D-{\rm div}D^{\rm sing})$, la
différentielle de $D$ induit une forme linéaire non nulle sur
l'espace tangent de $\bar g_{v}$. Il s'ensuit que la fraction d'éléments
à réduction $\bar g_v$
que l'on
obtient en déformant un $g_v$ donné  est $1/q_v$ de sorte que l'ensemble des
$g_{v}\in G^{\rm reg}(\Oc_{v})-G^{\rm tvl}(\Oc_{v})$ ayant la
réduction $\bar g_{v}$ est de mesure $q_v^{-\dim(G)-1}$. En plus,
les points de $G$ annulés par le discriminant forment un diviseur de
$G$. Puisque le diviseur $\mathbf c^{-1}({\rm div}D-{\rm div}D^{\rm sing})$ est de
codimension un, le complémentaire de $G^{\rm tvl}(\Oc_{v})$ dans
$G'(\kappa_{v})$ a aussi une mesure $O(q_{v}^{-2})$ lorsque $q_{v}$
tend vers l'infini.
Le lemme s'en déduit.
\end{proof}

Observons que d'après sa définition $G^{\rm tvl}$ est contenu dans 
$G^{\rm rs}(F_v)$, donc dans l'image réciproque $\mathbf c^{-1}(\mathfrak A^{\rm rs})(F_v)$.
Considérons la résolution simultanée de Gro\-then\-dieck-Springer
$\pi:\tilde G\to G$ dont la restriction à
$G^{\rm rs}$ est un morphisme fini $\pi^{\rm rs}:\tilde G^{\rm
rs}\to G^{\rm rs}$. Rattachés à un point de $\tilde G^{\rm rs}$ il y a un tore maximal
et un sous-groupe de Borel qui le contient. Il y a évidemment une identification
canonique de l'ensemble des sous-groupes de Borel qui contient un tore maximal donné
avec le groupe de Weyl.
L'image directe $\pi_{*}^{\rs} \Ql$ est donc un faisceau
$\ell$-adique muni d'une action de $W$.
Nous pouvons introduire
le faisceau
\begin{equation}\label{ngo}
(\pi_{*}^{\rs} \Ql \otimes X^{*}(T))^{W},
\end{equation}
mais nous avons besoin d'un objet plus délicat.

Il faut utiliser des résultats de l'article
\cite{N} démontrés sous l'hypothèse que la caractéristique
résiduelle ne divise pas l'ordre du groupe de Weyl, une hypothèse que
nous admettons pour les fins de cet article lacunaire. L'application de $g$ sur la partie semi-simple
de sa dé\-com\-po\-sition comme produit d'un élément semi-simple et d'un élément unipotent, donne
une application $G\to T/W$. Il y a aussi une application $\tilde G\to T$. Elles
donnent
\begin{equation*}
\begin{CD}
\tilde G @>>> T\\
@VVV @VVV\\
G @>>> T/W
\end{CD}
\end{equation*} 
Ce diagramme donne par restriction un diagramme cartésien
\begin{equation*}
\begin{CD}
\tilde G^{\rm reg} @>>> T\\
@VVV @VVV\\
G^{\rm reg} @>>> T/W
\end{CD}
\end{equation*}
L'argument
est le suivant. Le morphisme $G^{\rm reg} \to T/W$ est lisse et après un
changement de base le morphisme $G^{\rm reg} \times_{T/W} T \to T$
est aussi lisse. Par conséquent, le produit fibré est aussi lisse et donc normal.
On a un morphisme $\tilde G^{\rm reg}$ vers le produit fibré qui
est d'une part fini, car on sait que $\tilde G^{\rm reg} \to G^{\rm reg}$ est fini,
mais aussi birationnel car le diagramme est clairement cartésien au-dessus
de $G^{\rm rs}$. La normalité implique maintenant
que le morphisme
$$
\tilde G^{\rm reg} \to G^{\rm reg}\times_{T/W} T
$$
est un isomorphisme. Il en résulte que le groupe $W$ agit sur $\tilde G^{\rm reg}$; cette
action ne s'étend pas sur $\tilde G$.

Il est alors permis d'introduire le faisceau
$$
\mathcal L= (\pi_{*}^{\reg} \Ql \otimes X^{*}(T))^{W}.
$$
La restriction \eqref{ngo} de $\mathcal L$ à l'ouvert $G^{\rm rs}$ est
un système
local dont la monodromie est donnée par le revêtement
étale galoisien géométriquement $\tilde G^{\rs} \to G^{\rs}$ et
l'action de $W$ habituelle sur $X^{*}(T)\otimes \Ql$. On a supposé
$G$ semi-simple. Par conséquent la représentation de $W$
sur $X^*(T)$ ne contient pas de représentation triviale
et $\mathcal L$ n'admet pas un faisceau constant comme sous-quotient.

Puisque $G$ est simplement connexe, les centralisateurs des
éléments réguliers semi-simples sont connexes et la restriction du schéma des
centralisateurs $I$ sur $G$ à l'ouvert $G^{\rm rs}$ est un tore
$I\times_G G^{\rs}$. Les caractères de ce tore forment un système
local sur $G^{\rs}$ qui après le changement de coefficients à $\Ql$,
est isomorphe à $\mathcal L$. Nous n'avons pas trouvé de références
pour cette affirmation mais l'analogue dans le cas de l'algèbre de Lie a
été démontré par Donagi et Gaisgory dans \cite{DG}. La proposition 2.4.7 de
\cite{N} est aussi une référence commode. La démonstration dans le cas
d'un groupe simplement connexe est tout à fait semblable à celle pour
l'algèbre de Lie.

Si $g_v:\Spec(\Oc_v) \to G$ est un trait transversal au diviseur
discriminant, l'image réciproque de $\Spec(\Oc_v)$ dans $\tilde G$ est
une réunion de traits, autrement dit un revêtement fini et normal de
$\Spec(\Oc_v)$. L'argument justifiant cette affirmation peut être
trouvé dans la démonstration de 4.7.3 de \cite{N}. Le lemme
suivant, qui décrit la fibre de $\mathcal L_{\bar g_v}$ à la réduction $\bar g_v$
de $g_v$, s'en déduit.

\begin{lemme}
	Soit $g_{v}\in G^{\rm tvl}(\Oc_{v})$. Notons $X^{*}(T_{g_{v}})$ le
module des caractères du tore $T_{g_{v}}$. Il est défini sur une extension
algébrique finie de $F_{v}$ et muni d'une action continue du groupe de
Galois de $F_{v}$.
La partie $(X^{*}(T_{g_{v}})\otimes
\Ql)^{I_{v}}$ invariante sous le groupe d'inertie  est alors isomorphe 
à la fibre $\mathcal L_{\bar g_{v}}$
comme $\Ql$-espaces vectoriel munis d'une action de Frobenius.
\end{lemme}

Grâce à ce lemme nous pouvons introduire au lieu de la fonction $f_{v}(g_{v},s)$
la fonction
$f'_{v}(g_{v},s)$ supportée par $G^{\rm reg}(\Oc_{v})$ qui associe à
$g_{v}\in G^{\rm reg}(\Oc_{v})$ la valeur
\begin{align*}
f'_{v}(g_{v},s)&=1-\tr({\rm Fr}_{v}, \mathcal L_{{\bar g}_{v}}) q_v^{-1} 
+\tr({\rm Fr}_{v},\mathcal L_{{\bar g}_{v}}) q_v^{-s}\\
&=
1-\tr({\rm Fr}_{v}, \mathcal L_{{\bar g}_{v}})\left(q_v^{-1}-{q_v^{-s}}\right).
\end{align*}
Il suffit de comparer
l'intégrale des deux fonctions sur l'ouvert $G^{\rm tvl}(\Oc_{v})$
car l'intégrale de l'une ou de l'autre sur son complément est $O(q_v^{-2})$. 
Sur $G^{\rm tvl}(\Oc_{v})$ leur différence est aussi majorée par $q_v^{-2}$.

Pour démontrer la proposition \ref{prodinf}, nous observons d'abord
que
$$
\int_{G^{\rm tvl}(\Oc_v)}dg_v=\int_{G(\Oc_v)}dg_v+O(q_v^{-2})=
q_v^{-\dim G}|G(\kappa_v)|+O(q_v^{-2})
$$
et que
\begin{equation}\label{estimate1}
q_v^{-\dim G}|G(\kappa_v)|=1+O(q_v^{-2}).
\end{equation}
Cette dernière égalité se trouve dans les articles de Ono
sur les nombres de Tamagawa. 

Puisque la codimension de $G^{\rm reg}$ est $3$,
pour terminer la démonstration de la proposition \ref{prodinf}, il suffit de vérifier
la majoration,
\begin{equation}\label{estimate2}
q_v^{-\dim G}\sum_{\bar g_{v} \in G^{\rm reg}(\kappa_{v})}\tr({\rm Fr}_{v},\mathcal L_{{\bar g}_{v}})=
O(q_v^{-1/2}).
\end{equation}
C'est une conséquence immédiate du lemme suivant.
\begin{lemme}
       Lorsque $q_{v} \rightarrow \infty$, on a la majoration
$$\sum_{\bar g_{v} \in G^{\rm reg}(\kappa_{v})} {\rm tr}({\rm Fr}_{v},
\mathcal L_{{\bar g}_{v}}) =O(q_{v}^{\dim(G)-1/2}).$$
\end{lemme}

D'après la formule des traces de Grothendieck-Lef\-schetz
pour le faisceau $\mathcal L$ sur le schéma $G^{\rm reg}$ 
restreint à  ${\rm Spec}_{\kappa_v}$, l'expression
ci-dessus est égale à
$$
\sum_{i=0}^{2\dim(G)} (-1)^{i}{\rm tr}({\rm Fr},
\rmH^i_c(G^{\rm reg} \otimes_{\Oc_v} \bar \kappa_v, \mathcal L)).
$$
Notons que les dimensions des ces groupes de cohomologie sont
in\-dé\-pen\-dantes de la place $v$ à l'exception d'un nombre fini d'entre elles.
D'après  le théorème principal de \cite{D}, les valeurs absolues des
valeurs propres de Frobenius  sur 
$\rmH^{i}_{c}(G^{\rm reg}\otimes_{\mathcal O_v}\bar\kappa_v, \mathcal L)$ 
sont majorées par $q_{v}^{i/2}$. Il suffit donc de démontrer
que $\rmH^{2\dim(G)}_{c}(G^{\rm reg}, \mathcal L)=0$. Mais cela résulte du
fait que $\mathcal L$ n'admet pas de système local trivial comme
sous-quotient.

Pour le groupe $SL(2)$ et une caractéristique résiduelle impaire
la démonstration de la proposition \ref{prodinf} peut se faire
d'une façon éĺément\-aire. Il y a trois types de tore sur $F_v$: déployé, non ramifié et ramifié.
Leurs contributions à l'intégrale de la fonction constante $1$
sont de la forme 
$$
a_1+a_2q_v^{-1}+O(q_v^{-2}),\qquad
b_1+b_2q_v^{-1}+O(q_v^{-2}), \qquad c_2q_v^{-1}+O(q_v^{-2}).
$$ 
On vérifie à la main que $a_1+b_1=1$
et $a_2+b_2+c_2=0$.
Leurs contributions au coefficient de $q_v^{-1}-q_v^{-s}$ sont respectivement de la forme
$a_3$, $b_3$ et $0$. On vérifie, encore à la main, que $a_3+b_3=O(q_v^{-1})$. C'est
ce deuxième calcul qui exige un traitement bien plus raffiné --- la formule
de Grothendieck-Lef\-schetz --- dans le cas général.

\end{document}